\documentclass{amsart}

\usepackage{amsmath,amssymb}
\usepackage{amsthm}
\usepackage[nobysame]{amsrefs}
\usepackage{indentfirst}
\usepackage{mathrsfs}
\usepackage{graphicx}
\usepackage{hyperref}
\usepackage{xcolor}
\usepackage{fourier}
\usepackage[margin=1.5in]{geometry}
\usepackage[norelsize,boxed,linesnumbered,vlined,algo2e]{algorithm2e}

\usepackage{booktabs}
\usepackage{multirow}

\theoremstyle{plain}

\theoremstyle{definition}

\theoremstyle{remark}

\DeclareMathOperator{\tr}{tr}

\DeclareMathOperator{\rank}{rank}

\DeclareMathOperator{\diag}{diag}

\newcommand{\TT}{\mathbb{T}}

\newcommand{\CC}{\mathbb{C}}
\newcommand{\Or}{\mathcal{O}}

\IfFileExists{mathabx.sty}%
  {\DeclareFontFamily{U}{mathx}{\hyphenchar\font45}%
   \DeclareFontShape{U}{mathx}{m}{n}{<->mathx10}{}%
   \DeclareSymbolFont{mathx}{U}{mathx}{m}{n}%
   \DeclareFontSubstitution{U}{mathx}{m}{n}%
   \DeclareMathAccent{\widebar}{0}{mathx}{"73}%
}{%
  \PackageWarning{mathabx}{%
    Package mathabx not available, therefore\MessageBreak substituting
    widebar with overline\MessageBreak }%
  \newcommand{\widebar}[1]{\overline{#1}}%
}

\newcommand{\mc}[1]{\mathcal{#1}}

\newcommand{\abs}[1]{\lvert#1\rvert}

\newcommand{\average}[1]{\left\langle#1\right\rangle}

\newcommand{\omm}{\text{omm}}
\newcommand{\cond}{\text{cond}}

\title{Preconditioning orbital minimization method for planewave discretization}

\author{Jianfeng Lu} \address{Department of Mathematics, Department of
  Physics, and Department of Chemistry, Duke University, Box 90320, Durham NC 27708, USA}
\email{jianfeng@math.duke.edu}

\author{Haizhao Yang}
\address{Department of Mathematics, Duke University, Box 90320, Durham NC 27708, USA}
\email{haizhao@math.duke.edu}

\date{\today} \thanks{This work is partially supported by the National
  Science Foundation under grants DMS-1312659, DMS-1454939, and
  ACI-1450280. H.Y. thanks the support of the AMS-Simons Travel
  Award. We thank Lexing Ying for helpful discussions.}

\begin{document}

\begin{abstract}
  We present an efficient preconditioner for the orbital minimization
  method when the Hamiltonian is discretized using planewaves (i.e.,
  pseudospectral method). This novel preconditioner is based on an
  approximate Fermi operator projection by pole expansion, combined
  with the sparsifying preconditioner to efficiently evaluate the pole
  expansion for a wide range of Hamiltonian operators. Numerical
  results validate the performance of the new preconditioner for the
  orbital minimization method, in particular, the iteration number is
  reduced to $\Or(1)$ and often only a few iterations are enough for
  convergence.
\end{abstract}

\keywords{Kohn-Sham density functional theory, orbital minimization method, 
preconditioning, Fermi operator projection, pole expansion, sparsifying preconditioner}

{\bf AMS subject classifications.}  65K10, 65F08, 65F15

\maketitle

\section{Introduction}

We consider the problem of finding the low-lying eigenspace of a
Hamiltonian matrix $H$ coming from the planewave discretization of a
self-adjoint Hamiltonian operator $-\frac{1}{2}\Delta + V$ for some
potential function $V$. Given an $n \times n$ Hamiltonian matrix $H$,
it is well-known that, the eigenspace associated to the first $N$
eigenvalues (non-degeneracy is assumed throughout the work) is given
by the trace minimization
\begin{equation}\label{eq:tracemin}
  E = \min_{X \in \CC^{n \times N}, \, X^{\ast} X = I_N} \frac{1}{2}\tr ( X^{\ast} H X ), 
\end{equation}
where $I_N$ is a $N \times N$ identity matrix, thus $X^{\ast} X = I_N$
is the orthonormality constraints of the columns of $X$. 
When a
conventional algorithm is used to solve for the eigenvalue problem, a
QR factorization in each iteration is required to impose
orthogonality. The computational cost of each QR is $\Or(nN^2)$ with a
large prefactor and communication cost is expensive in high
performance computing. This creates an obstacle to minimize
\eqref{eq:tracemin} if the number of iteration is large. Much effort
has been devoted to reducing the communication cost, see e.g., 
\cites{Demmel:15,Martinsson:15} and references therein. 

It turns out that it is in fact possible to remove the orthogonality
constraint. In the context of linear scaling algorithms for
electronic structure calculations, the orbital minimization method
(OMM) was proposed in \cites{MauriGalliCar:93, MauriGalli:94,
  Ordejon:93, Ordejon:95} to circumvent the orthonormality constraint
in \eqref{eq:tracemin}. Instead, we search for the eigenspace by an
unconstrained minimization
\begin{equation}\label{eq:omm}
  E = \min_{X \in \CC^{n \times N}} E_{\omm}(X) = \min_{X \in \CC^{n \times N}}\frac{1}{2} \tr \bigl( ( 2 I_N - X^{\ast} X ) ( X^{\ast} H X ) \bigr).
\end{equation}
To see where \eqref{eq:omm} comes from, notice that we may reformulate \eqref{eq:tracemin} to relax the orthogonality constraint as 
\begin{equation}
  E = \min_{X \in \CC^{n \times N},\, \rank X = N}\frac{1}{2}\tr \bigl( (X^{\ast} X)^{-1} ( X^{\ast} H X) \bigr), 
\end{equation}
so that $X$ is no longer constrained. The factor
$(2 I_N - X^{\ast} X)$ in \eqref{eq:omm} can then be seen as a Neumann
series expansion of
$(X^{\ast} X)^{-1} = \bigl( I_N - (I_N - X^{\ast} X) \bigr)^{-1}$
truncated to the second term.

Somewhat surprisingly, for negative definite $H$ (all the eigenvalues
are strictly less than $0$), the minimizer of $E_{\omm}$ spans the
same subspace as the minimizer to the original problem
\eqref{eq:tracemin}, thus the desired eigenspace \cite{MauriGalliCar:93,MauriGalli:94}. Note that
the assumption that $H$ being negative definite can be made without
loss of generality, as we may just shift the diagonal of $H$ by a
constant, and this shift will not change the eigenspace. 

The OMM was originally proposed for linear scaling calculations (the
total computational cost is linearly proportional to $N$, the number
of electrons) for sparse $H$, combined with truncating $X$ to keep
only $\Or(1)$ entries per column. However, even without adopting the
linear scaling truncation, whose error can be hard to control, in the
context of cubic scaling implementation, the OMM algorithm still has
advantage over direct eigensolver in terms of scalability in parallel
implementations \cite{Corsetti:14}.  The OMM has the potential to be a
competitive strategy for finding the eigenspace.

Since the OMM transforms the eigenvalue problem into an unconstrained
minimization, the nonlinear conjugated gradient method is usually
employed to minimize the energy \eqref{eq:omm}. The efficiency of the
OMM algorithm thus depends on the optimization scheme, which in turn
crucially depends on the preconditioner.

In this work, we will consider $H$ coming from a Fourier
pseudospectral discretization of the Hamiltonian operator. The
pseudospectral method typically requires minimal degree of freedoms
for a given accuracy among standard discretizations and is also easy
to implement using the fast Fourier transform (FFT), and hence widely
used in the field physics and engineering literature. In the field of
electronic structure calculation, it is known as the planewave
discretization, which is arguably the most popular discretization
scheme up-to-date.

While the importance of preconditioning for planewave discretization
is well known and has been long recognized \cite{Teter:89}; the preconditioner
has been mostly limited to the type of shifted inverse Laplacian (more
details later). A natural question is whether more efficient
preconditioner can be designed.

We revisit the issue of preconditioning for planewave
discretization. We propose a new preconditioner using an approximate
Fermi operator projection based on the pole expansion (see e.g.,
\cites{Baroni:92,Lin:09,Polizzi:09,Goedecker:95, Goedecker:99}).  Once
constructed, the new preconditioner can be applied efficiently and
reduces the number of iteration of the OMM to only a few. The
resolvents involved in the pole expansion are solved iteratively using
GMRES \cite{SaadSchultz:86} combined with the recently developed
sparsifying preconditioner \cites{Ying:spspc, Ying:spspd,
  LuYing:nonSch}. Since only an approximate Fermi operator is
required, the construction and application of preconditioner become
quite efficient. Hence, the overall preconditioned OMM requires much
less computational time compared to existing OMM algorithms for
planewave discretizations.

This paper is organized as follows. In Section \ref{sec:POMM},
existing preconditioned OMMs are introduced and analyzed to motivate
the design of the new preconditioner. In Section \ref{sec:aglo}, we
construct the new preconditioner base on the approximate Fermi
operator projection and the sparsifying preconditioner. In Section
\ref{sec:numerics}, some numerical examples are provided to
demonstrate the efficiency of the new preconditioner. We conclude with
a discussion and some future work in Section \ref{sec:conclusion}.

\section{Preconditioned orbital minimization method}
\label{sec:POMM}

\subsection{Orbital minimization method}
The orbital minimization method (OMM) has become a popular tool in
electronic structure calculations for solving the Kohn-Sham eigenvalue
problem.  One of the major advantage is that unlike conventional
methods for eigenvalue calculations, no orthogonalization is required
in the iteration. The OMM minimizes the functional $E_{\omm}$, recalled here:
\begin{equation*}
  E = \min_{X \in \CC^{n \times N}} E_{\omm}(X) = \min_{X \in \CC^{n \times N}} \frac{1}{2} \tr \bigl( ( 2 I_N - X^{\ast} X ) ( X^{\ast} H X ) \bigr).
\end{equation*}
The nonlinear conjugate gradient method is usually applied for this
unconstrained minimization problem.  When the minimization problem
becomes ill-conditioned (for example, when the spectral gap, the
difference between the $N$-th and $(N+1)$-th eigenvalues, is small;
further discussed below), it requires a large number of iterations to
achieve convergence. Hence, an efficient preconditioner is crucial to
this minimization problem.

Let us consider preconditioning the gradient of $E_{\omm}$ in the
general framework of nonlinear conjugate gradient methods. Denote the
gradient as
\begin{equation}
  \mc{G}(X) := \frac{\delta E_{\omm}(X)}{\delta X^{\ast}} = 2 H X - X (X^{\ast} H X) - HX (X^{\ast} X).
\end{equation} 
The preconditioned nonlinear conjugate gradient method can be summarized as follows.

\begin{algorithm2e}[H]
\label{alg:pcg}
\caption{Preconditioned nonlinear conjugate gradient method.}
Initialize: Pick initial guess $X_1$ and fix a preconditioner $\mc{P}$; 

Set $D_1 = -\mc{P} \mc{G}(X_1)$ and perform a line search in the direction of
$D_1$ and update $X_2$; 

Set $m = 2$;

\While{not converged}{ Calculate the preconditioned gradient
  direction: $G_m = - \mc{P} \mc{G}(X_m)$;

  Compute $\beta_m$ according to the Polak-Ribi\`ere formula (other
  choices are equally possible)
  \begin{equation*}
    \beta_m = \frac{ G_m^{\ast} ( G_m - G_{m-1})}{G_{m-1}^{\ast} G_{m-1}};
  \end{equation*}

  Update the conjugate direction $D_m = G_m + \beta_m G_{m-1}$; 

  Perform a line search in the direction of $D_m$ and update $X_{m+1}$; 
  
  Set $m = m + 1$. 
}
\end{algorithm2e}

The OMM method is highly efficient in each iteration, as involves only
matrix-matrix multiplications that take $\Or(nN^2 + Nn\log n)$
operations, where $n\log n$ comes from the FFT in the application of
$H$ (recall that a pseudospectral discretization is used).  Hence, a
good preconditioner for the OMM should be efficient to construct and
simple to apply, so that it does not increase much the cost of each
iteration; at the same time, we hope the preconditioner can reduce the
number of iterations significantly.

\subsection{Existing preconditioners}
We shall first recall existing preconditioners for planewave
discretization in general, which also apply for the OMM.  The
conventional preconditioning employed in the electronic structure
calculation is in a form of inverse shifted Laplacian:
\begin{equation*}
  \mc{P} = P \otimes I_N, \quad \text{where,} \quad P =  (I - \tau^{-1} \Delta)^{-1}. 
\end{equation*}
Here $-\Delta$ stands for the discretization of the Laplacian operator
(kinetic energy operator in the Hamiltonian), and $\tau$ is a parameter
setting the scale for the kinetic energy preconditioning. Note that in
the planewave discretization, $P$ is a diagonal matrix in the
$k$-space
\begin{equation}
\label{eqn:Lap}
  P_{kk'} = \delta_{kk'} (1 + \tau^{-1} \abs{k}^2)^{-1}, 
\end{equation}
and thus
application of the preconditioner has minimal computational cost. 

In practice, the preconditioner might take a more complicated
expression (see e.g., the review article \cite{Payne:92} and
references therein)
\begin{equation}
\label{eqn:TPA}
  P_{kk'} = \delta_{kk'} \frac{27 + 18 s + 12 s^2 + 8s^3}{27 + 18s + 12s^2 + 8s^3 + 16 s^4}
\end{equation}
with $s = \abs{k}^2 / \tau$ and $\tau$ a scaling parameter. This
preconditioner was first proposed in \cite{Teter:89} and is now known
as the TPA preconditioner. It has a similar asymptotic behavior to
\eqref{eqn:Lap} around $\abs{k} = 0$ and as $\abs{k} \to \infty$, but
is found to be more efficient than \eqref{eqn:Lap} empirically.

A recent article \cite{Zhou:15} generalized the TPA
preconditioner using polynomials with higher degrees than those in
\eqref{eqn:TPA}. More specifically, this new preconditioner, denoted
as the generalized TPA preconditioner (gTPA), is based on a polynomial
of degree $t$
\[
p_t(s):=c_0+c_1s+c_2s^2+\cdots +c_ts^t,
\]
and takes the form
\begin{equation}
\label{eqn:gTPA}
 P_{kk'} = \delta_{kk'} \frac{p_t(s)}{p_t(s)+c_{t+1}s^{t+1}},
\end{equation}
where $s = \abs{k}^2 / \tau$ as before. The polynomial $p_t$ is
constructed such that all the derivatives of $g_t(s) =
\frac{p_t(s)}{p_t(s)+c_{t+1}s^{t+1}}$ up to order $t$ at $s=0$ vanish, meaning that the gTPA preconditioner is close
to an identify operator in the low frequency regime and the width of
this region depends on $t$. Note that, the TPA preconditioner is a special
case of the gTPA preconditioners for $t=3$ and the shifted Laplacian
preconditioner in \eqref{eqn:Lap} is a special case of $t=0$. From
now on, $P_t$ is used to denote the gTPA preconditioner corresponding
to $t$.
In the pseudospectral method, applying all the preconditioners above to a search direction only takes $\Or(N n\log n)$ operations since FFT is applied to $N$ columns in the search direction $ \mc{G}(X_m)$.

\begin{figure}
\begin{center}
   \begin{tabular}{c}
        \includegraphics[height=2.5in]{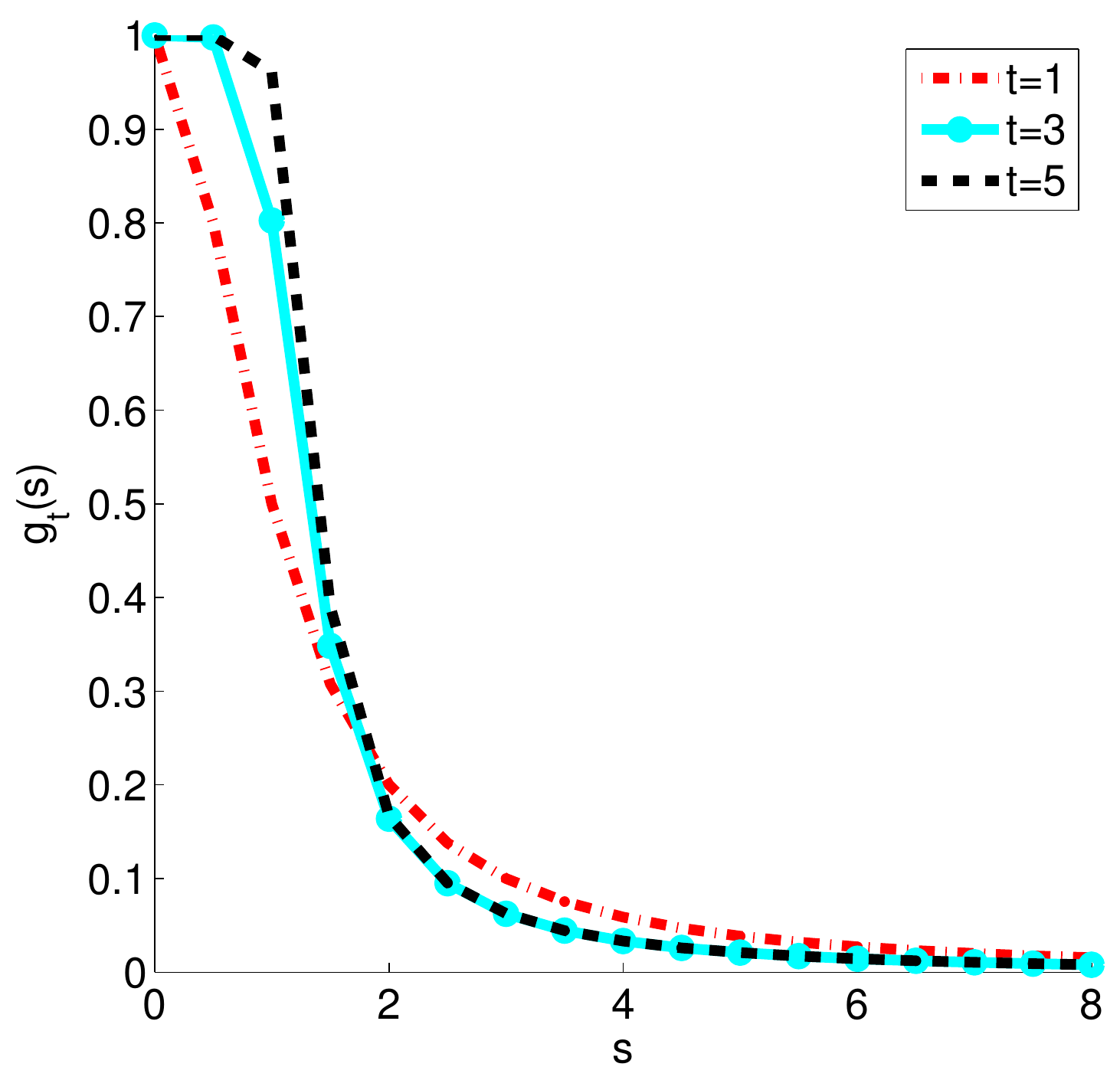} 
\end{tabular}
\end{center}
\caption{The function $g_t(x)$ in the gTPA preconditioner for $t=1$, $3$, and $5$.}
\label{fig:p}
\end{figure}

The physical intuition of all these preconditioners is that: 1) the preconditioner should keep the low-frequency components in the search direction unchanged; 2) $g_t$ should have an asymptotic behavior like $\frac{1}{s}$ so that it behaves as an inverse Laplacian for the high-frequency components. Therefore, they act as a low-pass filter that is essentially similar to the inverse shifted Laplacian in \ref{eqn:Lap}. For the gTPA, as $t$ increases, the low-frequency region is expanded while the high frequency components are further damped (see Figure \ref{fig:p} for an illustration). 

As the behavior of these preconditioners is close to inverse shifted
Laplacian, the success essentially lies in the assumption that the
potential part of the Hamiltonian operator does not change the
spectrum too much. When this assumption is not valid, these
preconditioner might not be able to effectively reduce the condition
number of the OMM. 

\subsection{Condition number of the preconditioned OMM}
To understand the performance of the existing preconditioners in the
context of OMM, we analyze here the condition number of the
preconditioned OMM. While the condition number of the OMM (without
preconditioning) was analyzed in \cite{PfrommerDemmelSimon:99}, for
completeness, we will also repeat the calculation here in a slightly
different and more compact way. 

Note that condition number of the TPA type preconditioning is
rather complicated (if possible) to obtain analytically. For the purpose of simplicity, we
consider a preconditioner of the type of inverse shifted Hamiltonian
\begin{equation}
\label{eqn:invH}
P = (H-\mu I)^{-1},
\end{equation}
which is more friendly for analysis (since it commutes with the
Hamiltonian operator) and shares the same spirit with the TPA type
preconditioners.  In fact, as it takes into account both the kinetic
and potential parts of the Hamiltonian, it captures better the
spectral behavior of the Hamiltonian operator.

The Hessian of $E_{\omm}$ is given by
\begin{equation}
  \mc{H}(X) Z = 2 H Z - Z (X^{\ast} H X) - X (Z^{\ast} H X) - X (X^{\ast} H Z) 
  - H Z (X^{\ast} X) - HX (Z^{\ast} X) - HX (X^{\ast} Z).
\end{equation}
Let $X_0$ be the minimizer (i.e., the solution of the eigenvalue problem) such that $H X_0 = X_0 \Lambda_0$, where $\Lambda_0$ is a diagonal matrix containing the eigenvalues. 
Evaluating the Hessian at the minimizer $X_0$ (orthonormalized such that $X_0^{\ast} X_0 = I_N$), we have
\begin{equation}
  \mc{H}(X_0) Z = H Z - Z \Lambda_0 - X_0 Z^{\ast} X_0 \Lambda_0  - X_0 \Lambda_0 Z^{\ast} X_0 - 2 X_0 \Lambda_0 X_0^{\ast} Z. 
\end{equation}
Note that the energy $E_{\omm}$ is invariant under any perturbation in the subspace of $X_0$, we may assume that the perturbed direction is perpendicular, i.e., $Z^{\ast} X_0 = 0$. For such directions, we get 
\begin{equation}
  \mc{H}(X_0) Z = H Z - Z \Lambda_0. 
\end{equation}
Therefore, the condition number of $\mc{H}$ is at least 
\begin{equation}
  \cond(\mc{H}(X_0)) \geq \frac{ \lambda_n - \lambda_1}{\lambda_{N+1} - \lambda_N},
\end{equation}
determined by the ratio
of the spectrum width of $H$ and the spectral gap.
In particular, this implies that if the gap between $\lambda_N$ and
$\lambda_{N+1}$ is small relatively to the whole spectrum width of
$H$, the condition number of the OMM minimization becomes large.

For the preconditioned gradient, we get
\begin{equation}
  \mc{P} \mc{G}(X) = 2 P H X - PX (X^{\ast} H X) - P H X (X^{\ast} X),  
\end{equation}
and thus
\begin{equation}
  \mc{P} \mc{H}(X_0) Z = P H Z - P Z \Lambda_0.
\end{equation}
The eigenvalues of the Hessian hence depend on $P (H -
\lambda_i I)$, $i = 1, \ldots, N$, and some of them take the form
\begin{equation}
\label{eqn:eigv}
  1 - \frac{ \lambda_i - \mu}{\lambda_j - \mu}, \qquad i \leq N < j. 
\end{equation}

To guarantee that the Hessian after preconditioning is positive
definite, we consider the family of preconditioner with $P = (H - \mu
I)^{-1}$ with $\mu < \lambda_{N+1}$. To estimate the condition number
after preconditioning, we separate into several cases according to the
position of $\mu$ in the spectrum of $H$.
\begin{itemize}
\item $\lambda_N< \mu<\lambda_{N+1}$, then the largest and the smallest eigenvalues of the form \eqref{eqn:eigv} after preconditioning are 
\[
\lambda_{\max} = 1-\frac{\lambda_1-\mu}{\lambda_{N+1}-\mu},
\]
\[
\lambda_{\min} = 1-\frac{\lambda_N - \mu}{\lambda_{n}-\mu},
\]
since $i\leq N$ and $N+1\leq j$.
Then the condition number is at least 
\begin{equation}
  \cond(\mc{P}\mc{H}(X_0)) \geq \frac{1 - \frac{ \lambda_1 - \mu}{\lambda_{N+1} - \mu}}{1 - \frac{\lambda_N - \mu}{\lambda_{n} - \mu}} = \frac{ \lambda_{N+1} - \lambda_1}{\lambda_{n} - \lambda_N} \frac{ \lambda_{n} - \mu }{\lambda_{N+1} - \mu}. 
\end{equation}
To minimize the condition number, we let $\mu\rightarrow \lambda_{N}$ and get 
\begin{equation}
  \cond(\mc{P}\mc{H}(X_0)) \geq  \frac{ \lambda_{N+1} - \lambda_1}{\lambda_{N+1} - \lambda_{N}}. 
\end{equation}
\item $\lambda_1< \mu<\lambda_{N}$, then
\[
\lambda_{\max} = 1-\frac{\lambda_1-\mu}{\lambda_{N+1}-\mu},
\]
\[
\lambda_{\min} = 1-\frac{\lambda_N - \mu}{\lambda_{N+1}-\mu},
\]
since $i\leq N$ and $N+1\leq j$.
Then the condition number is at least
\begin{equation}
  \cond(\mc{P}\mc{H}(X_0)) \geq  \frac{1 - \frac{ \lambda_1 - \mu}{\lambda_{N+1} - \mu}}{1 - \frac{\lambda_N - \mu}{\lambda_{N+1} - \mu}} = \frac{ \lambda_{N+1} - \lambda_1}{\lambda_{N+1} - \lambda_N}. 
\end{equation}
\item $\mu\leq \lambda_1$, then
\begin{equation}
  \cond(\mc{P}\mc{H}(X_0)) \geq  \frac{1 - \frac{ \lambda_1 - \mu}{\lambda_n - \mu}}{1 - \frac{\lambda_N - \mu}{\lambda_{N+1} - \mu}} = \frac{ \lambda_n - \lambda_1}{\lambda_{N+1} - \lambda_N} \frac{\lambda_{N+1} - \mu}{\lambda_n - \mu}.
\end{equation}
To minimize the condition number, we should take $\mu$ close to $\lambda_1$ in this case, which leads to 
\begin{equation}
  \cond(\mc{P}\mc{H}(X_0)) \geq \frac{ \lambda_{N+1} - \lambda_1}{\lambda_{N+1} - \lambda_N}. 
\end{equation}
\end{itemize}

Therefore, we arrive at the conclusion that the lower bound of the
conditioner we can achieve by using a preconditioner of type $P = (H -
\mu I)^{-1}$ with $\mu <\lambda_{N+1}$ is \[\frac{ \lambda_{N+1} -
  \lambda_1}{\lambda_{N+1} - \lambda_N}.\] Hence, the condition number
after preconditioning will still be large if the spectral gap between
$\lambda_N$ and $\lambda_{N+1}$ is small.

The above conclusion can be understood by investigating the spectral
meaning of an inverse shifted Hamiltonian of the form $(H-\mu
I)^{-1}$, we see that this type of preconditioners prefers the
eigenspace with eigenvalues close to the shift $\mu$. This would limit
the search direction in this spectral window and hence the OMM cannot
converge to the right solution quickly, or even converges to an 
undisired solution. This motivates us to propose a
new preconditioner based on an approximate Fermi operator projection
that restricts the search direction in the full target
eigenspace. This algorithm is presented in the next section.

\section{Algorithm description}
\label{sec:aglo}

The preconditioner we propose is based on the idea of using an
approximate Fermi operator projection, as described in
Section~\ref{sec:contour}. The preconditioner has the form of a linear
combination of shifted inverses of the Hamiltonian; so that to accelerate the
construction and application, the sparsifying preconditioner
\cite{Ying:spspd} is used to iteratively solve the linear equation
corresponds to the shifted inverse. The detail is given in
Section~\ref{sec:sparse}. Last, it is possible to precompute and store
the preconditioner in a data sparse format using rank revealing QR
decomposition, as will be discussed in Section~\ref{sec:precompute}. 

\subsection{Fermi operator projection}\label{sec:contour}

In quantum mechanics, given an effective one-particle Hamiltonian $H$, the zero-temperature single-particle density matrix $\Pi$ of the system is given by the Fermi operator via a Green's function expansion
\begin{equation}\label{eqn:FOP}
\Pi=\frac{1}{2\pi i}\int_{\Gamma}(z I-H)^{-1}dz,
\end{equation}
where $\Gamma$ is a contour in the complex plane containing the $N$
eigenvalues below the Fermi level. Note that by Cauchy's integral
formula, $\Pi$ defined in \eqref{eqn:FOP} is exactly the projection to
the eigenspace corresponds to the first $N$ eigenvalues.  The contour
integral representation of the Fermi operator has been a useful tool
in electronic structure calculation, for example in linear scaling
methods (see e.g., the review \cite{Goedecker:99}). More recently, the contour integral
formulation was used in the method of PEXSI \cites{Lin:09,LinLuYingCarE:09,Lin:11} for a
fast algorithm to obtain the density for sparse Hamiltonian matrices,
combined with the selected inversion algorithm.  It is also used in
FEAST \cite{Polizzi:09} as a general eigenvalue solver for sparse
Hamiltonian matrices. In \cite{Damle:14}, the contour integral was
used for multishift problems.  Our idea in this work is to explore the
approximate Fermi operator projection as an effective preconditioner.

To use $\Pi$ as an preconditioner, similar to what we did in
\cite{Lin:09}, we employ the efficient quadrature rule proposed in
\cite{Hale:08} to discretize the contour integration \eqref{eqn:FOP}
as
\begin{equation}
\label{eqn:dFOP}
\Pi \approx \sum_{j=0}^p w_j (H-z_j I)^{-1}
\end{equation}
where $z_j$ is the $j$-th quadrature node on the contour and $w_j$ is
its corresponding weight. The above formula is also called pole
expansion, as the quadrature points can be viewed as poles in the
resolvents $(H-z_j I)^{-1}$ in the expansion. The approximation error
of the quadrature rule is bounded by \cite{Lin:09}
\begin{equation*}
  e^{-cp \bigl(  \log \frac{\lambda_n - \lambda_1}{\lambda_{N+1} - \lambda_N} + 3 \bigr)^{-1}}.
\end{equation*} 
The convergence of the quadrature rule is exponentially fast in $p$
and hence only a small number of poles are needed. The details of the
choice of $w_j$ and $z_j$ can be found in \cite{Lin:09} and we will
not repeat here.

We remark that if the projection is calculated by inverting $(H-z_j
I)^{-1}$ exactly (or with a very small error tolerance), we get
already the projection $\Pi$ onto the low-lying eigenspace, i.e., the
density matrix. Of course, direct inversion of $(H-z_j I)^{-1}$, for
pseudospectral discretization, is equally expensive, since $H$ is a
dense matrix. This is different from the situation of PEXSI or FEAST
where sparse Hamiltonian matrices are considered.  An alternative
approach to construct a precise projection is to use iterative matrix
solvers to apply $(H-z_j I)^{-1}$. If a basis of the eigenspace is
desired, it can be obtained by acting $\Pi$ on a few vectors, as in
the Fermi operator projection method \cites{Baroni:92, Goedecker:99}, or more
systematically, by using a low-rank factorization via a randomized SVD
\cite{Halko:11}. However, to obtain high accuracy, this idea is still
expensive since the matrix $(H - z_j I)$ is usually ill-conditioned
and the iteration number might be large.

The key observation is that, as we plan to apply the projection as a
preconditioner, we in fact do not really need the exact projection, as
long as we can get a good approximation. Thus, it suffices to use an
iterative scheme to solve $(H-z_j)^{-1}$ acting on some vector with a
relatively large error tolerance. We aim to achieve a balance, such
that the approximate projection is accurate enough such that the OMM
converges in a few iterations, but also rough enough such that the
construction is cheap.

More specifically, in the nonlinear conjugate gradient method for the
OMM, the preconditioned gradient direction at the $m$-th step is
computed by
\begin{equation}
\label{eqn:direct}
G_m = - \mc{P} \mc{G}(X_m) =  -\sum_{j=0}^p w_j Y_{m,j},
\end{equation}
where $Y_{m,j}\in \mathbb{R}^{n\times N}$ is a rough solution of the linear system
\begin{equation}
\label{eqn:lst}
(H-z_j I)Y_{m,j}=\mc{G}(X_m),
\end{equation}
for $j=1$, $\dots$, $p$.  In the pseudospectral method, the matrix
$H-z_j I$ can be applied efficiently via the
FFT in $\Or(n\log n)$ operations. Hence, the GMRES
\cite{SaadSchultz:86} with a large tolerance and a small number of
maximum iterations is able to provide rough solutions to the linear
systems in \eqref{eqn:lst} quickly. Denote the number of iterations in
GMRES by $n_g$, the total complexity for finding a rough solution is
$\Or(p n_g N n\log n)$, which is of the same order as the complexity
of the regular TPA and gTPA preconditioners. To further accelerate the
GMRES iteration, the search direction $\mc{G}(X_{m})$ is used as the
initial guess for the GMRES, since it is almost in the eigenspace of
interest if $m$ is large. As we shall see later, the new
preconditioner is highly efficient and the preconditioned OMM will
converge in only a few iterations.

\subsection{Sparsifying preconditioner}\label{sec:sparse}

In the previous subsection, an approximate Fermi operator projection is introduced as a preconditioner for the OMM. As we shall see later in the numerical section, the preconditioned OMM converges quickly in a few iterations. Hence, the difficulty of solving an eigenvalue problem with the OMM has been replaced by the difficulty of constructing the projection efficiently, the most expensive part of which is solving linear systems with multiple right-hand sides. Since these linear systems are  solved roughly in an iterative scheme like the GMRES, the preconditioning problem for the OMM is transferred to the preconditioning problem for the GMRES.

In the case when the Hamiltonian operator $H=-\frac{1}{2}\Delta+V-\mu$ behaves like a kinetic energy operator, i.e., the potential energy operator is negligible, a conventional preconditioner for the linear system 
\begin{equation}
\label{eqn:lstsp}
 (H-\mu I)u = b
\end{equation}
is 
\[
P = (-\frac{1}{2}\Delta -s)^{-1},
\]
where $s$ is a proper shift. However, when the contribution of the potential energy operator $V$ is large, e.g., in highly indefinite systems on periodic structures, an inverse shifted Laplacian is no longer an efficient preconditioner, e.g., the GMRES might take a very large number of iterations to converge or it even diverges. To
overcome this difficulty, we will adopt the recently proposed sparsifying
preconditioner for solving the linear systems in the construction of the approximate Fermi operator projection.

Let us focus on the numerical solution of \eqref{eqn:lstsp}. 
We assume that the computation domain is the periodic unit box
$[0,1]^d$ and discretize the problem using the Fourier pseudospectral
method. 
With abuse of notations, the discretized problem of
\eqref{eqn:lstsp} takes the form

\begin{equation}\label{eq:Disc}
  \bigl(-\frac{1}{2}\Delta + V - \mu\bigr) u = b.
\end{equation}
We will briefly recall the key idea of the sparsifying preconditioner
for solving this type of equations. More details can be found in
\cite{Ying:spspd} (see also \cites{Ying:spspc, LuYing:nonSch}).

Denote $l = \average{V}$ the spatial average of $V$. We assume
without loss of generality that $-\frac{1}{2}\Delta + l - \mu$ is invertible
on $\TT^d$ with periodic boundary condition, otherwise, we will use a
slight perturbation of $\mu$ instead and put the difference into
$V$. This allows us to rewrite \eqref{eqn:lstsp} trivially as
\begin{equation}
  \bigl(-\frac{1}{2}\Delta + (l-\mu)  + (V-l) \bigr) u = b
\end{equation}
Applying the Green's function of the constant part $G =
\bigl(-\frac{1}{2}\Delta + (l-\mu) \bigr)^{-1}$ via the FFT
to both sides of \eqref{eqn:lstsp}, we have an equivalent
linear system
\begin{equation}\label{eq:GLu}
  \bigl(I + G (V - l) \bigr) u = G b. 
\end{equation}

The main idea of the sparsifying preconditioner is to multiply a
particular sparse matrix $Q$ (to be defined later) on the both hand
sides of \eqref{eq:GLu}:
\begin{equation}\label{eq:QG}
  (Q + QG (V - l) ) u = Q (Gb), 
\end{equation}
so to make the matrix on the left hand in \eqref{eq:QG} becomes
sparse. To see how this is possible, let $S =\{(j,a(j)), j\in J\}$ be
the support of the sparse matrix $Q$, i.e., for each point $j$, the
row $Q(j,:)$ is supported in a local neighborhood $a(j)$.  If for each
point $j$, the $Q$ is constructed such that 
\begin{equation*}
  (QG)(j,a(j)^c)=Q(j,a(j))G(a(j),a(j)^c) \approx 0,
\end{equation*}
then the product $QG$ is also essentially supported on $S$. The sum
$Q+QG (V-l)$ is essentially supported in $S$ as well, since $(V-l)$ is
a diagonal matrix. The above requirement can indeed be achieved since
$G$, after all, is the Green's function of a local operator; thus $Q$
is a discrete approximation of the differential operator. More details can be found in \cite{Ying:spspd}.


Restricting $Q+QG (V-l)$ to $S$ by thresholding
other values to zero, we define a sparse matrix $P$
\begin{equation}\label{eq:P}
  P_{ij}=
  \begin{cases}
    \bigl(Q+QG (V-l)\bigr)_{ij}, & (i,j)\in S,\\
    0, &  (i,j)\not\in S.
  \end{cases}
\end{equation}
Since $P\approx Q+QG (V-l)$, we have the approximate equation
\begin{equation}
  P u \approx Q(Gb). 
\end{equation}
The sparsifying preconditioner computes an approximate solution
$\tilde{u}$ by solving
\[
P \tilde{u} = Q(Gb).
\]
Since $P$ is sparse, the above equation can be solved by sparse direct
solvers such as the nested dissection algorithm \cite{George:73}. 
The solution $\tilde{u}=P^{-1} QG b$
can be used as a preconditioner for the standard iterative algorithms
such as GMRES \cite{SaadSchultz:86} for the solution of
\eqref{eq:Disc}.


In the construction of the sparsifying preconditioner, the most
expensive part is to build the nested dissection factorization for
$P$, the complexity of which scales as $\Or(n^{1.5})$ in two
dimensions and $\Or(n^2)$ in three dimensions. The application of the
sparsifying preconditioner is very efficient and takes only $\Or(n\log
n)$ operations in two dimensions and $\Or(n^{4/3})$ operations in
three dimensions. Since the dominant cost of the OMM with
pseudospectral discretization is $\Or(nN^2)$, where $N$ is
proportional to $n$, the construction and the application of the
sparsifying preconditioner is relatively cheap for large-scale
problems. As we shall see in the numerical section, even if $n$ is as
small as $576$, the preconditioned OMM with the sparsifying
preconditioner is still more efficient than the one with a regular
preconditioner.

\subsection{Precomputing the preconditioner}\label{sec:precompute}

When the problem size $n$ is not sufficiently large, the prefactor $pn_g$ in the application complexity of the approximate projection might be too large. It is better to use a randomized low-rank factorization method to construct a matrix representation $UU^*$ of the preconditioner such that 
\[
UU^*b\approx \sum_{j=0}^p w_j (H-z_j I)^{-1}b,
\]
for an arbitrary vector $b\in \mathbb{R}^{n}$, where $U\in
\mathbb{C}^{n\times N}$ is a unitary matrix (and hence the notation
$U$). Therefore, $UU^{\ast}$ is a data-sparse way to store the matrix
preconditioner. Adapting the randomized SVD method \cite{Halko:11}, a
fast algorithm for constructing the matrix $U$ is given below.

\begin{algorithm2e}[H]
\label{alg:QR}
\caption{Randomized low-rank factorization for the Fermi operator projection.}
Given the Hamiltonian $H$ and the range of the eigenvalues, construct the quadrature nodes and weights $\{z_j,w;j\}_{1\leq j\leq p}$ for the contour integration in \ref{eqn:FOP};

Construct a Gaussian random matrix $B\in \mathbb{R}^{n\times (N+k)}$, where $k$ is a non-negative constant integer;

Solve the linear systems $(H-z_j I)Y_j = B$ for each $j=1$, $\dots$, $p$ using the GMRES with the right-hand side as the initial guess.

Compute $Y =  \sum_{j=0}^p w_j Y_j$;

Apply a rank-revealing QR factorization $[U,R]=\mathsf{qr}(Y)$;

Update $U$ by selecting the first $N$ columns: $U = U(:,1:N)$;

The approximate Fermi operator projection is given by $UU^*$.
\end{algorithm2e}

The operation complexity of Algorithm \ref{alg:QR} is $\Or(p n_g N
n\log n+nN^2)$, where $nN^2$ comes from the QR factorization. When $N$
is very large and $nN^2$ dominates the construction complexity of the
preconditioner, it is better to apply the projection in
\eqref{eqn:direct} directly without the QR factorization, while for
smaller scale problems, the pre-computation accelerates the overall
calculation. 

We remark that the scalability of the rank-revealing QR factorization in
high performance computing has been an active research direction, see
for example the communication-reduced QR factorization recently
proposed in \cite{Martinsson:15} and implemented in Elemental
\cite{Poulson:2013}. Due to the large prefactor of the complexity of
the communication-reduced QR factorization, applying it for many times
is still quite expensive. This makes Algorithm \ref{alg:QR} better
than those with QR factorization in each iteration, e.g., in
conventional ways to solve for \ref{eq:tracemin}. Parallelism of other
steps is straightforward: 1) each linear system at each pole $z_j$ can
be solved simultaneously; 2) parallel GMRES routines have been
standard in high performance computing; 3) applying this approximate
Fermi operator projection is simple matrix-matrix multiplication with
complexity $\Or(nN^2)$.

\section{Numerical examples}\label{sec:numerics}

This section presents numerical results to support the efficiency of the  proposed algorithm. Numerical results were obtain in MATLAB on a Linux computer with CPU speed at 3.5GHz. The number of poles in the contour integration was 30. The GMRES algorithm in the construction of the approximate Fermi operator projection was used with a relative tolerance equal to $10^{-5}$, a restart number equal to 15, and a maximum iteration number equal to 5 (i.e., the total maximum iterations allowed in the GMRES is $75$). In the randomized algorithm for constructing a matrix representation $UU^*$ of the approximate Fermi operator projection, a random Gaussian matrix $B$ of size $n\times N$ was used. In the OMM, the convergence tolerance was $10^{-13}$ and the maximum iteration number is $4000$. 

Numerical experiments corresponding to three different kinds of
Hamiltonian operators in two dimensions are presented: 1) the
potential energy is much weaker than the kinetic energy, i.e., the
kinetic energy operator dominates the Hamiltonian operator; 2) the
potential energy is prominent, contains one defect, 
and the Hamiltonian operator behaves different from 
the kinetic operator; 3) the potential energy is more 
prominent and contains more defects, so that the Hamiltonian operator
behaves further different from the kinetic operator.  These
experiments demonstrate that the proposed preconditioner has better
performance than the regular preconditioners in a wide range of
circumstances.

In our numerical tests, the contour construction is based on the exact
spectral information of the $H$. The starting guess for the conjugate
gradient procedure in the OMM is generated by diagonalizing the
Hamiltonian $H$, selecting the low-lying eigenpairs, and adding
Gaussian random noise with a distribution $\mathcal{N}(0,0.1M^2)$ to
the eigenvectors, where $M$ is the maximum magnitude of the ground
truth eigenvectors. Since the proposed preconditioner is also an approximate 
spectral projector onto the desired eigen subspace, we applied it to filter 
the initial guess before the preconditioned conjugate 
gradient scheme. In practice, spectral information can be estimated
in the previous step of the self-consistent field iteration in
electronic structure calculation.

\begin{figure}
\begin{center}
   \begin{tabular}{cc}
        \includegraphics[height=2in]{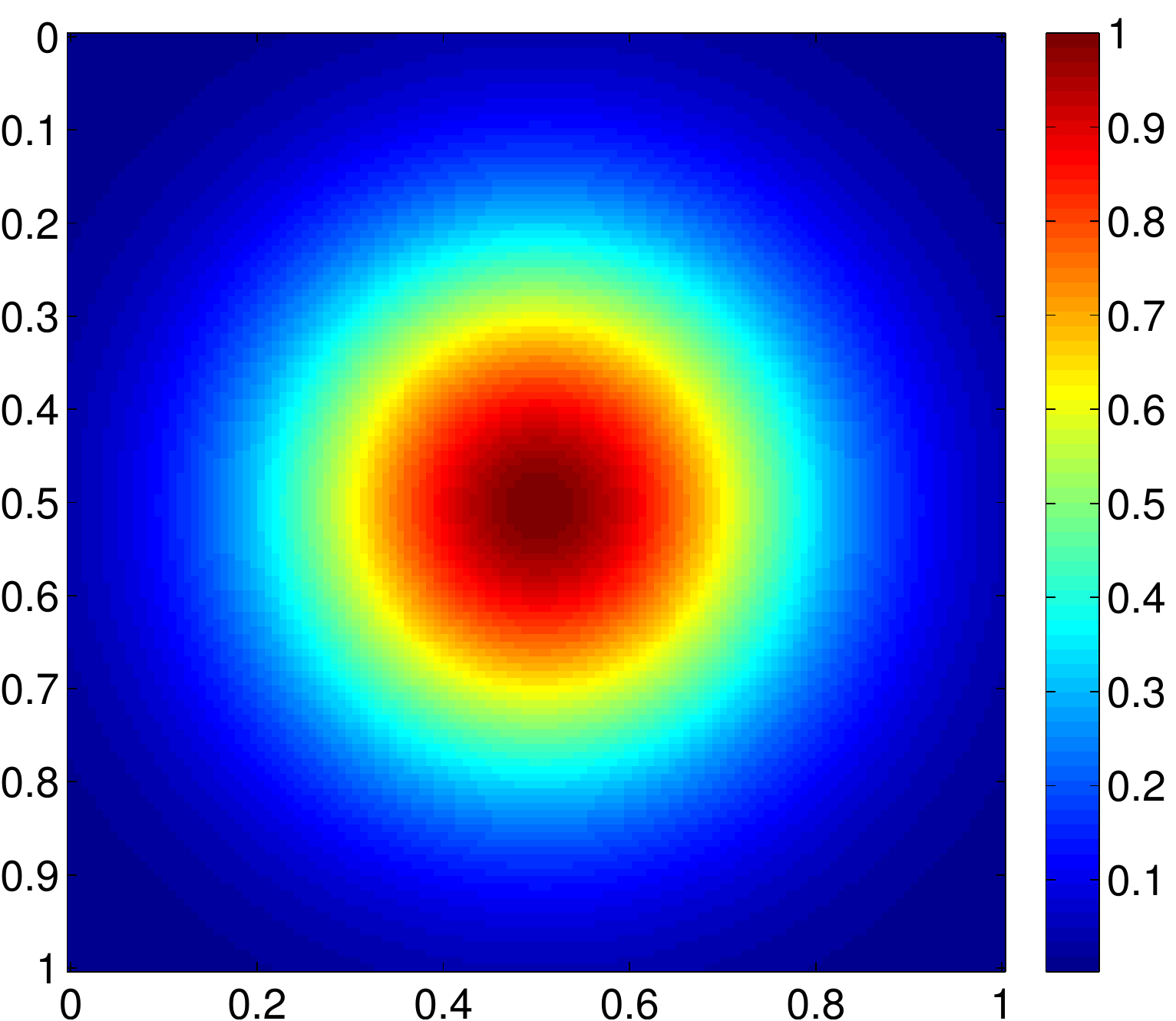}& \includegraphics[height=2in]{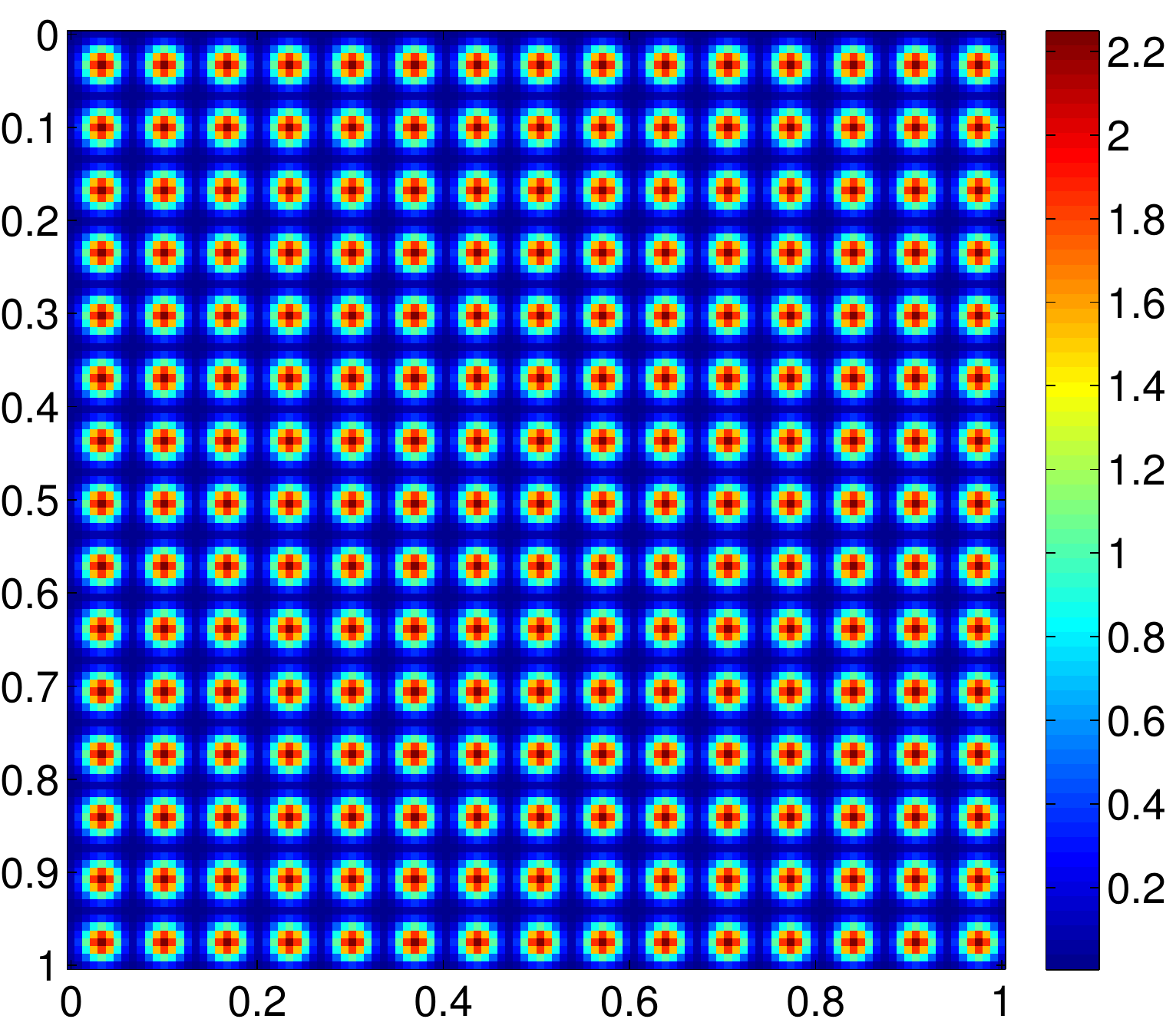}  
  \end{tabular}
\end{center}
\caption{Left: the periodic function $V_0(\mathbf{r})$ is a Gaussian well on the unit square $[0,1)^2$. Right: the potential energy operator $ \ell^2 V(\ell\mathbf{x})$ on the unit square $[0,1)^2$ after scaling with $\ell=15$ in Test $1$.}
\label{fig:pot1}
\end{figure} 

\subsection{Test $1$}\label{sec:test1}
In the first test, the Hamiltonian matrix $H$ is a discrete representation of the Hamiltonian operator in two dimensions
\begin{equation}
\label{eqn:ex0}
\left( -\frac{\Delta}{2} + V(\mathbf{r})\right) \phi_j(\mathbf{r})=\epsilon_j\phi_j(\mathbf{r}),\qquad \mathbf{r}\in \ell \mathbb{T}^2 :=[0,\ell)^2,
\end{equation}
with a periodic boundary condition, where $V(\mathbf{r})$ is the potential field, $\epsilon_j$ is the orbital energy of the corresponding Kohn-Sham orbital, $\phi_j(\mathbf{r})$. It is convenient to rescale the system to the unit square via the transformation $ \ell \mathbf{x} =\mathbf{r}$: 
\begin{equation}
\label{eqn:ex1}
\left( -\frac{\Delta}{2} + \ell^2 V(\ell\mathbf{x})\right) \phi_j(\mathbf{x})=\epsilon_j \ell^2 \phi_j(\mathbf{x}),\qquad \mathbf{x}\in \mathbb{T}^2 :=[0,1)^2,
\end{equation}
and discretize the new system with the pseudospectral method. Define 
\begin{equation*}
  J = \bigl\{(j_1, j_2) \mid 0 \leq j_1, j_2 < \sqrt{n} \bigr\},
\end{equation*}
\begin{equation*}
  K = \bigl\{ (k_1, k_2) \mid -\sqrt{n}/2 \leq k_1, k_2 < \sqrt{n}/2 \bigr\}.
\end{equation*}
Suppose $F$ and $F^{-1}$ are the discrete Fourier and inverse Fourier
transforms. After discretization, the corresponding eigenvalue problem of Equation \eqref{eqn:ex1} becomes 
\begin{equation}
\label{eqn:ex1mat}
HX=X\Lambda
\end{equation}
in numerical linear algebra, where $H=- \frac{\Delta}{2}+V$,  
\begin{equation*}
  - \frac{\Delta}{2}= F^{-1} \diag(2 \pi^2 \abs{k}^2)_{k \in K} F, 
\end{equation*}
and 
\begin{equation*}
  V = \ell^2 \diag\left(V\left(\frac{\ell j}{\sqrt{n}}\right)\right)_{j \in J},
\end{equation*}
$\Lambda$ is a diagonal matrix containing eigenvalues, and $X$ contains the eigenvectors.

\begin{figure}
\begin{center}
   \begin{tabular}{cc}
        \includegraphics[height=2in]{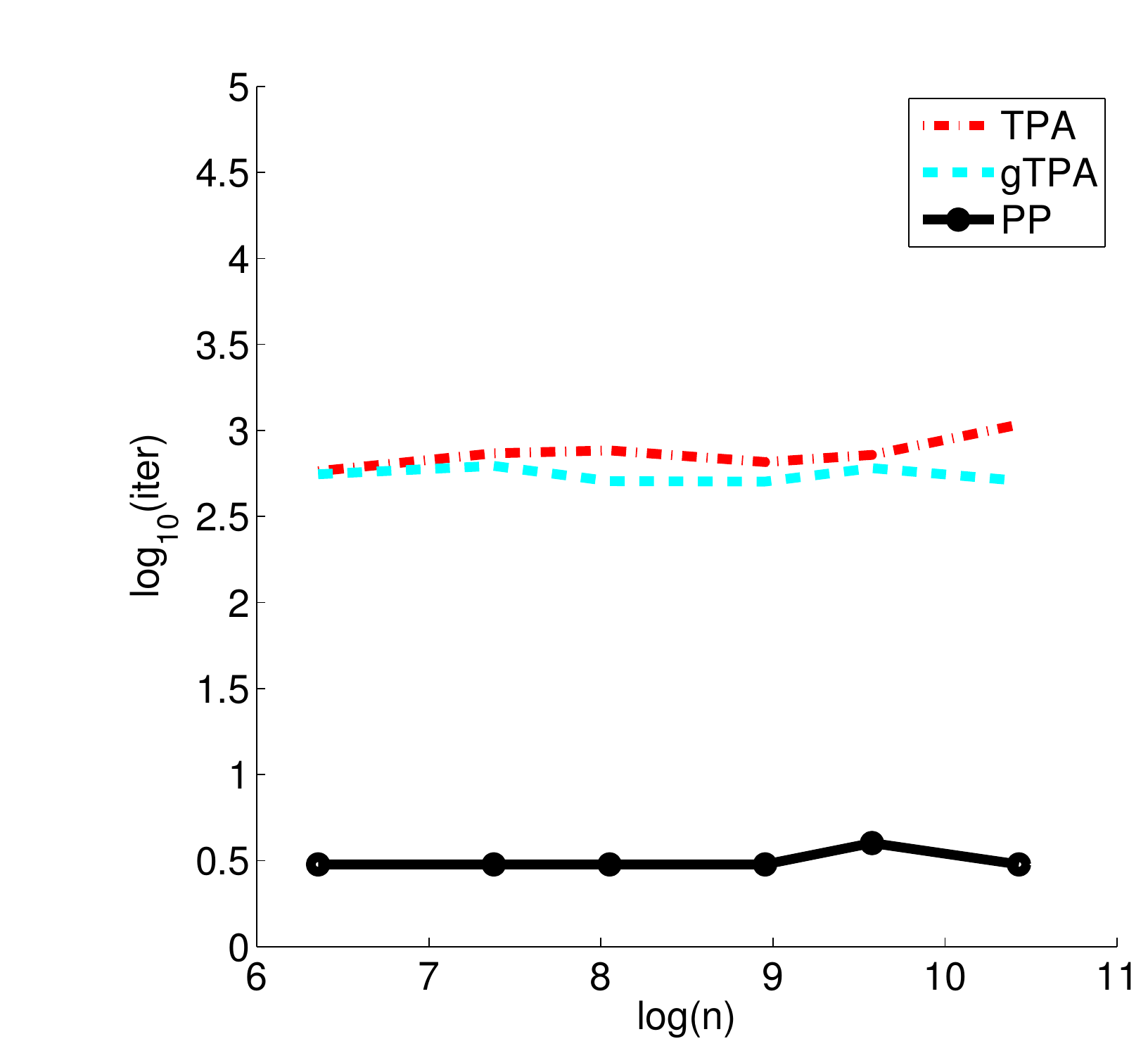} &\includegraphics[height=2in]{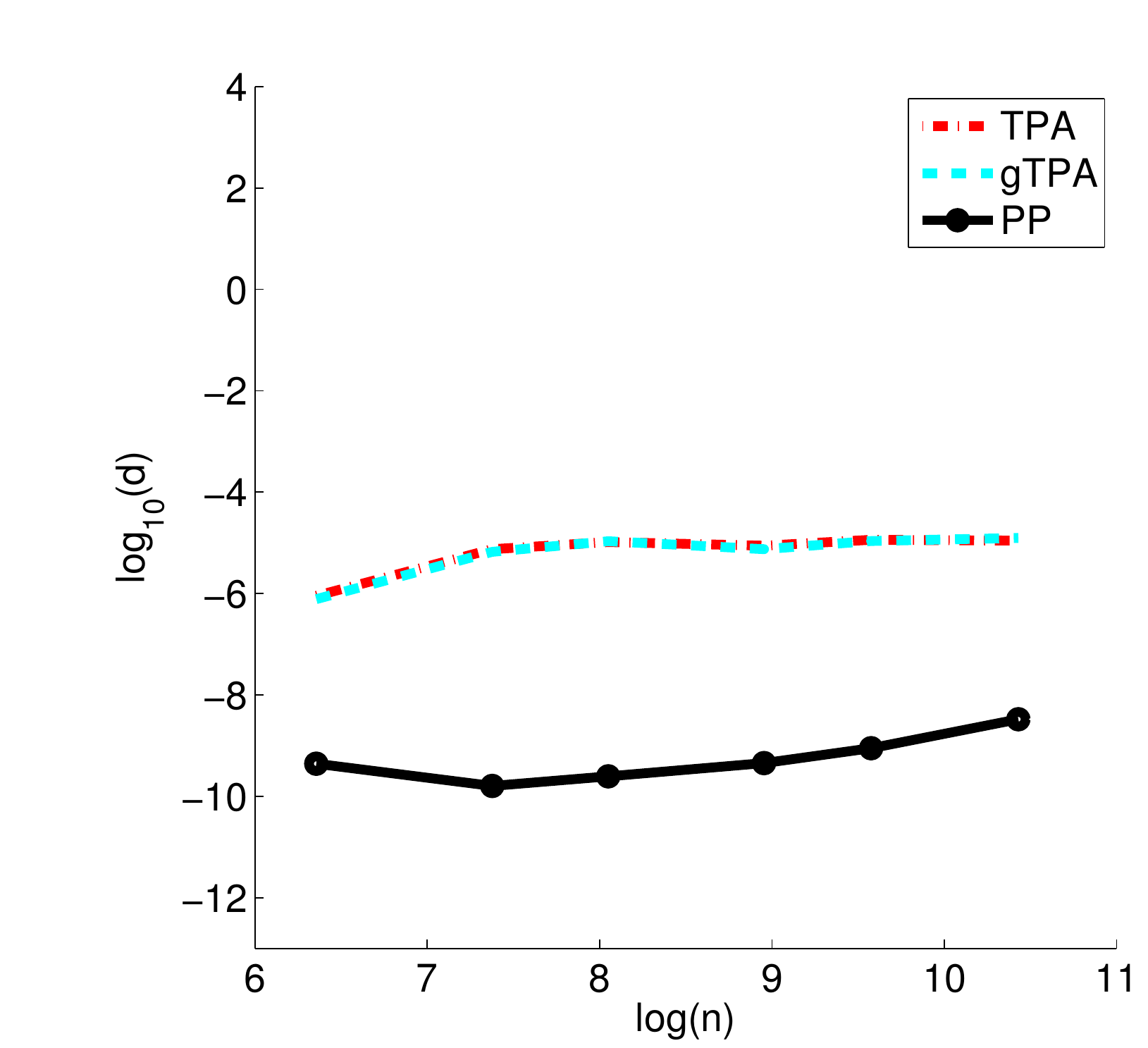}  \\
        (a) & (b)  \\
         \includegraphics[height=2in]{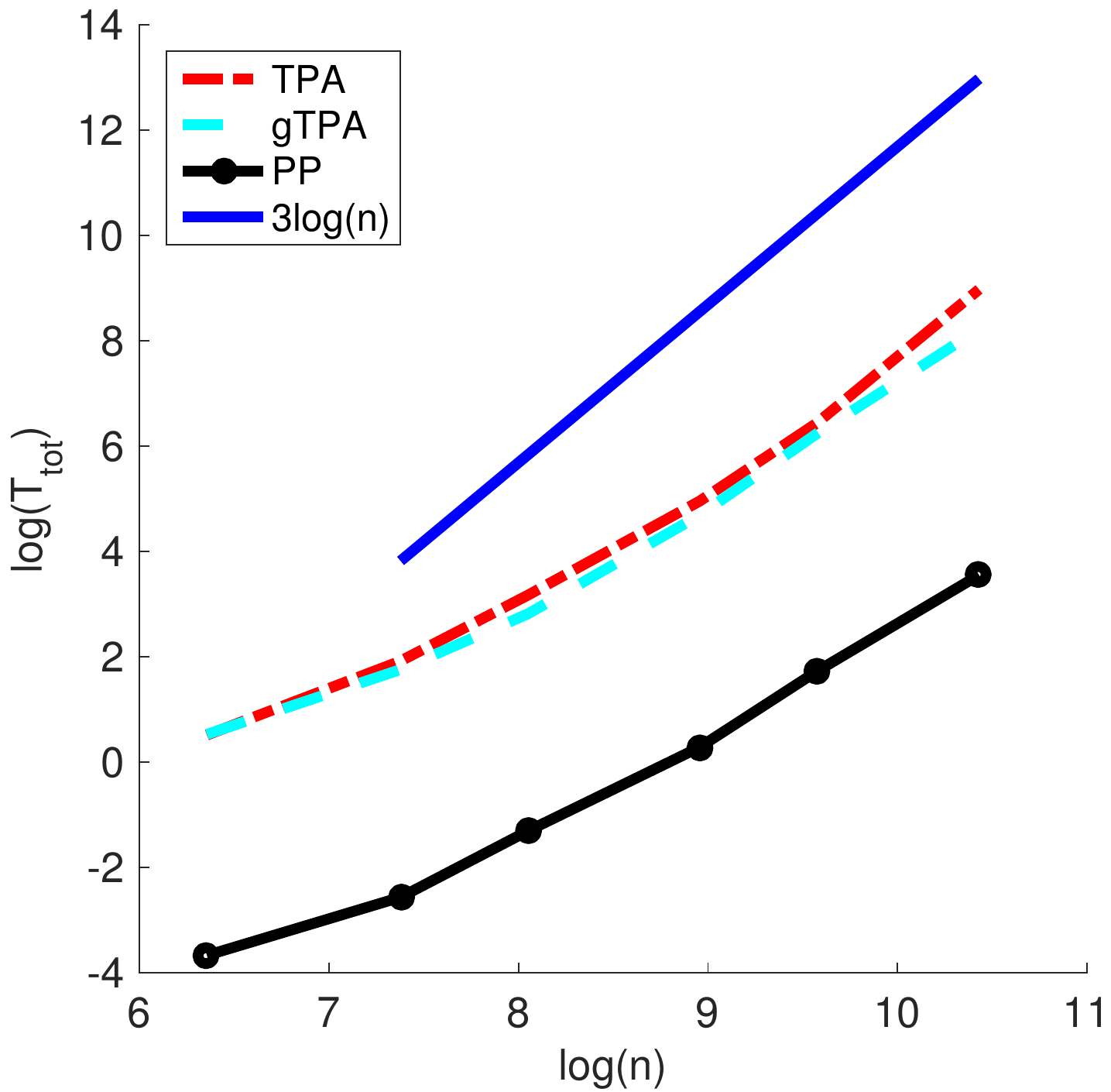}& \includegraphics[height=2in]{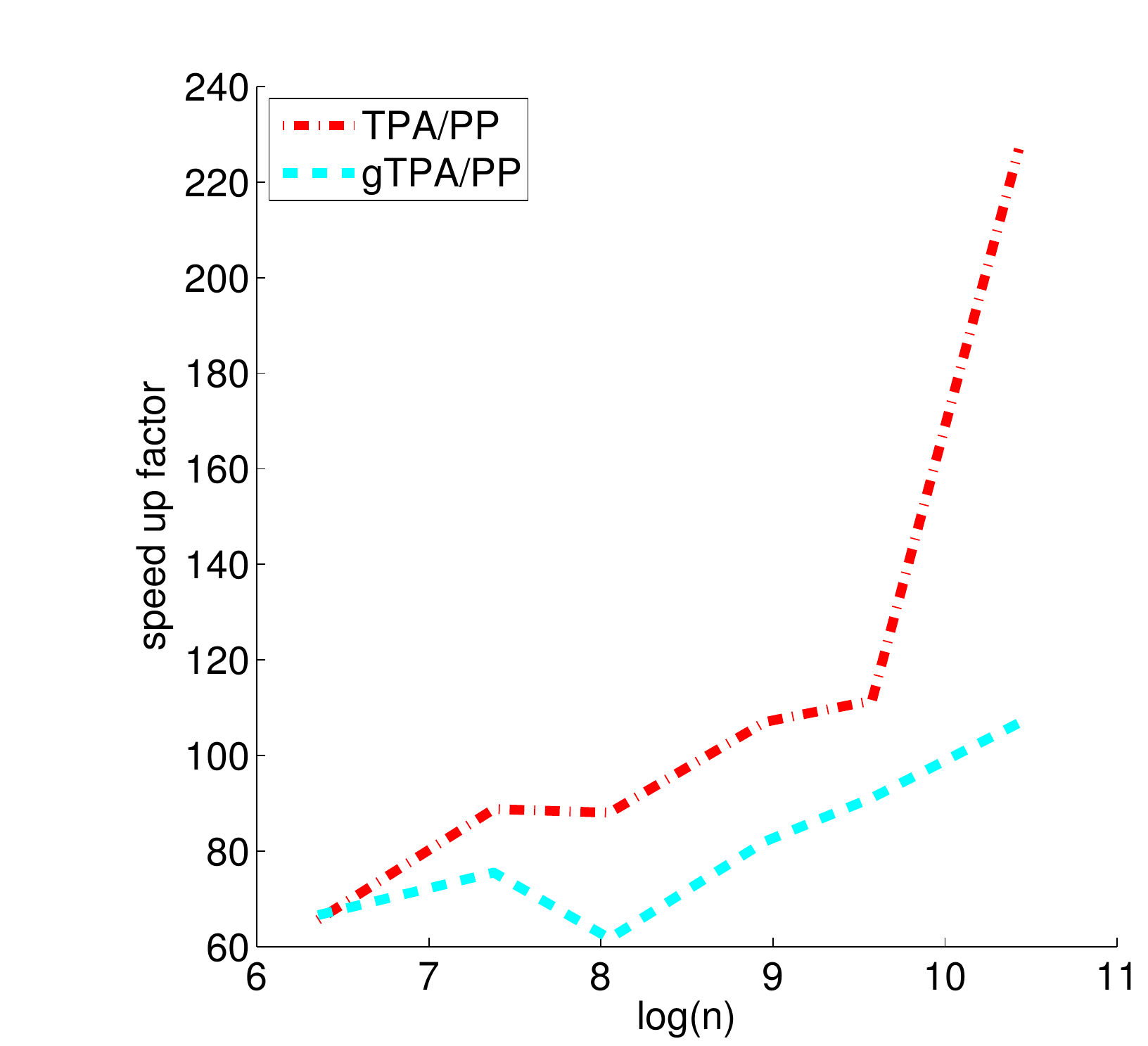} \\
         (c) & (d)
\end{tabular}
\end{center}
\caption{Numerical results in Test $1$ when $V(\mathbf{r})=0.01V_0(\mathbf{r})$. (a) the number of iterations in the preconditioned OMM. (b) the measurement $d(X,X_0)$. (c) the total running time $T_{\text{tot}}$. (d) the speedup factor of the PP method compared with the TPA and gTPA method.}
\label{fig:set1}
\end{figure}

Let $V_0(\mathbf{r})$ be a Gaussian well on the unit square $[0,1)^2$
(see Figure \ref{fig:pot1} left panel) and extend it periodically with
period $1$ in both dimensions. In the first test, the potential term
$V(\mathbf{r})$ (see Figure \ref{fig:pot1} right panel) in \eqref{eqn:ex1}
is chosen to be $0.01V_0(\mathbf{r})$ so that the kinetic energy
operator $ -\frac{\Delta}{2}$ dominates the Hamiltonian
operator. Hence, the spectrum of $H$ is very close to the one of $
-\frac{\Delta}{2}$. The preconditioned OMM is applied to solve the
eigenvalue problem in \eqref{eqn:ex1mat}. The numerical performance
of the empirical preconditioner in \eqref{eqn:TPA} (denoted as TPA)
and its generalization in \eqref{eqn:gTPA} when $t=5$ (denoted as
gTPA), and the new preconditioner via an approximate Fermi operator
projection (denoted as SPP or PP with or without the sparsifying
preconditioner, respectively), is compared and summarized in Table
\ref{tab:1} and Figure \ref{fig:set1}. The scaling parameter $\tau$ in
the TPA and gTPA preconditioner takes the value $\max_{j}\sum_{k}
\frac{1}{2}|k|^2(x^{(j)}(k))^2$, where $x^{(j)}$ is the $j$th column
of ground true eigenspace $FX_0$ and $k$ is the two-dimensional index
in $K$. For each Hamiltonian $H$ with a fixed $\ell$, $5$ experiments
were repeated with different random initial guesses. In this summary,
some notations are introduced and recalled as follows:
\begin{itemize}
\item $n$ is the dimension of $H$;
\item $\ell$ is the number of cells in the domain of $V(\mathbf{r})$ and each cell $[0,1)^2$ contains a grid of size $8\times 8$; $N=\ell$;
\item $\mathsf{cond}$ is the condition number of the OMM;
\item $\mathsf{iter}$ is the number of iterations in the preconditioned OMM;
\item $T_{\text{st}}$ is the setup time of the preconditioner; for PP and SPP, $T_{\text{st}}$ is the setup time per pole.
\item $T_{\text{omm}}$ is the running time of the preconditioned OMM;
 \item $T_{\text{tot}} = T_{\text{st}}+T_{\text{omm}}$;
\item $d$ measures the distance between the column space of the ground truth eigenvectors $X_0$ and the one of the estimated eigenvectors $X$ by the preconditioned OMM:
\[
d(X,X_0) =\frac{ \bigl\lVert X(X^*X)^{-1}X^* - X_0(X_0^*X_0)^{-1}X_0^* \bigr\rVert}{\bigl\lVert  X_0(X_0^*X_0)^{-1}X_0^* \bigr\rVert}, 
\]
where the norm here is the entrywise $\ell^\infty$ norm. 
\end{itemize}

\begin{table}[htp]
\centering
\begin{tabular}{rccccccc}
  \toprule
  & $(\ell,n)$ & $\mathsf{cond}$ & $\mathsf{iter}$ & $T_{\text{st}}(sec)$
  &  $T_{\text{omm}}(sec)$   & $T_{\text{tot}}(sec)$ & $d$ \\
\toprule
TPA & (3,576) & 1.4e+02 & 5.8e+02 & 2.271e-03 & 1.669e+00 & 1.672e+00 & 9.0e-07 \\ 
gTPA & (3,576) & 1.4e+02 & 5.6e+02 & 2.067e-03 & 1.695e+00 & 1.697e+00 & 7.8e-07 \\ 
PP & (3,576) & 1.4e+02 & 3.0e+00 & 1.486e-02 & 1.063e-02 & 2.549e-02 & 4.4e-10 \\ 
\toprule
TPA & (5,1600) & 8.0e+02 & 7.4e+02 & 4.380e-04 & 6.821e+00 & 6.822e+00 & 7.5e-06 \\ 
gTPA & (5,1600) & 8.0e+02 & 6.2e+02 & 5.310e-04 & 5.801e+00 & 5.802e+00 & 6.6e-06 \\ 
PP & (5,1600) & 8.0e+02 & 3.0e+00 & 4.204e-02 & 3.483e-02 & 7.686e-02 & 1.6e-10 \\ 
\toprule
TPA & (7,3136) & 1.6e+03 & 7.7e+02 & 1.596e-03 & 2.387e+01 & 2.388e+01 & 1.0e-05 \\ 
gTPA & (7,3136) & 1.6e+03 & 5.1e+02 & 2.333e-03 & 1.672e+01 & 1.673e+01 & 1.1e-05 \\ 
PP & (7,3136) & 1.6e+03 & 3.0e+00 & 1.464e-01 & 1.249e-01 & 2.713e-01 & 2.5e-10 \\ 
\toprule
TPA & (11,7744) & 1.3e+03 & 6.6e+02 & 4.097e-03 & 1.417e+02 & 1.417e+02 & 8.8e-06 \\ 
gTPA & (11,7744) & 1.3e+03 & 5.0e+02 & 4.640e-03 & 1.086e+02 & 1.086e+02 & 7.5e-06 \\ 
PP & (11,7744) & 1.3e+03 & 3.0e+00 & 6.151e-01 & 7.096e-01 & 1.325e+00 & 4.5e-10 \\ 
\toprule
TPA & (15,14400) & 7.2e+03 & 7.2e+02 & 1.296e-02 & 6.263e+02 & 6.263e+02 & 1.1e-05 \\ 
gTPA & (15,14400) & 7.2e+03 & 6.0e+02 & 1.245e-02 & 5.120e+02 & 5.121e+02 & 1.1e-05 \\ 
PP & (15,14400) & 7.2e+03 & 4.0e+00 & 1.622e+00 & 3.997e+00 & 5.619e+00 & 8.9e-10 \\ 
\toprule
TPA & (23,33856) & 1.7e+04 & 1.1e+03 & 6.347e-02 & 7.920e+03 & 7.920e+03 & 1.1e-05 \\ 
gTPA & (23,33856) & 1.7e+04 & 5.1e+02 & 6.488e-02 & 3.727e+03 & 3.727e+03 & 1.2e-05 \\ 
PP & (23,33856) & 1.7e+04 & 3.0e+00 & 9.595e+00 & 2.531e+01 & 3.490e+01 & 3.3e-09 \\

\bottomrule
\end{tabular}
\caption{Numerical results in Test $1$ when $V(\mathbf{r})=0.01V_0(\mathbf{r})$.}
\label{tab:1}
\end{table}

The preconditioner TPA and gTPA essentially assume that the potential $V(\mathbf{r})$ is negligible and they can be implemented efficiently using the fast Fourier transform. Hence, they would be efficient in the numerical examples in Test $1$ and this is supported by the results below: $T_{\text{st}}$ is negligible; the preconditioned OMM converges in a few hundred iterations; the OMM returns an eigenspace close to the ground truth, i.e., the measurement $d$ is almost $0$. 

Even in this case which favors the conventional preconditioners, the
new preconditioner via an approximate Fermi operator projection (PP)
is much more efficient. Even though the approximation accuracy of the
Fermi operator projection is low, since the tolerance of the GMRES in
computing the contour integration is large, this rough projection as a
preconditioner is able to reduce the iteration number of the OMM to
$3$ or $4$ in all examples. The total running time of PP, including
the setup time of the preconditioner and the running time of the OMM,
is far less than those of TPA and gTPA. The speedup factor is
significant and increases as the problem gets larger in most cases as
show in Figure~\ref{fig:set1}. The measurement $d$ is much smaller
than those in the TPA and gTPA methods, which means that the solution
of the OMM is more accurate than the TPA and gTPA methods.  We remark
that the solutions for linear systems at different poles can be
straightforwardly parallelized, here we only use the setup time per
pole in the comparison and focus on the time spent on the iterative
procedures of GMRES and OMM that can not be parallelized. Since the
GMRES with a regular preconditioner is already rather efficient and the
iteration converges less than $5$ iterations, the sparsifying
preconditioner is not applied in Test $1$.

\subsection{Test $2$}
          
\begin{figure}
\begin{center}
   \begin{tabular}{cc}   
            \includegraphics[height=2in]{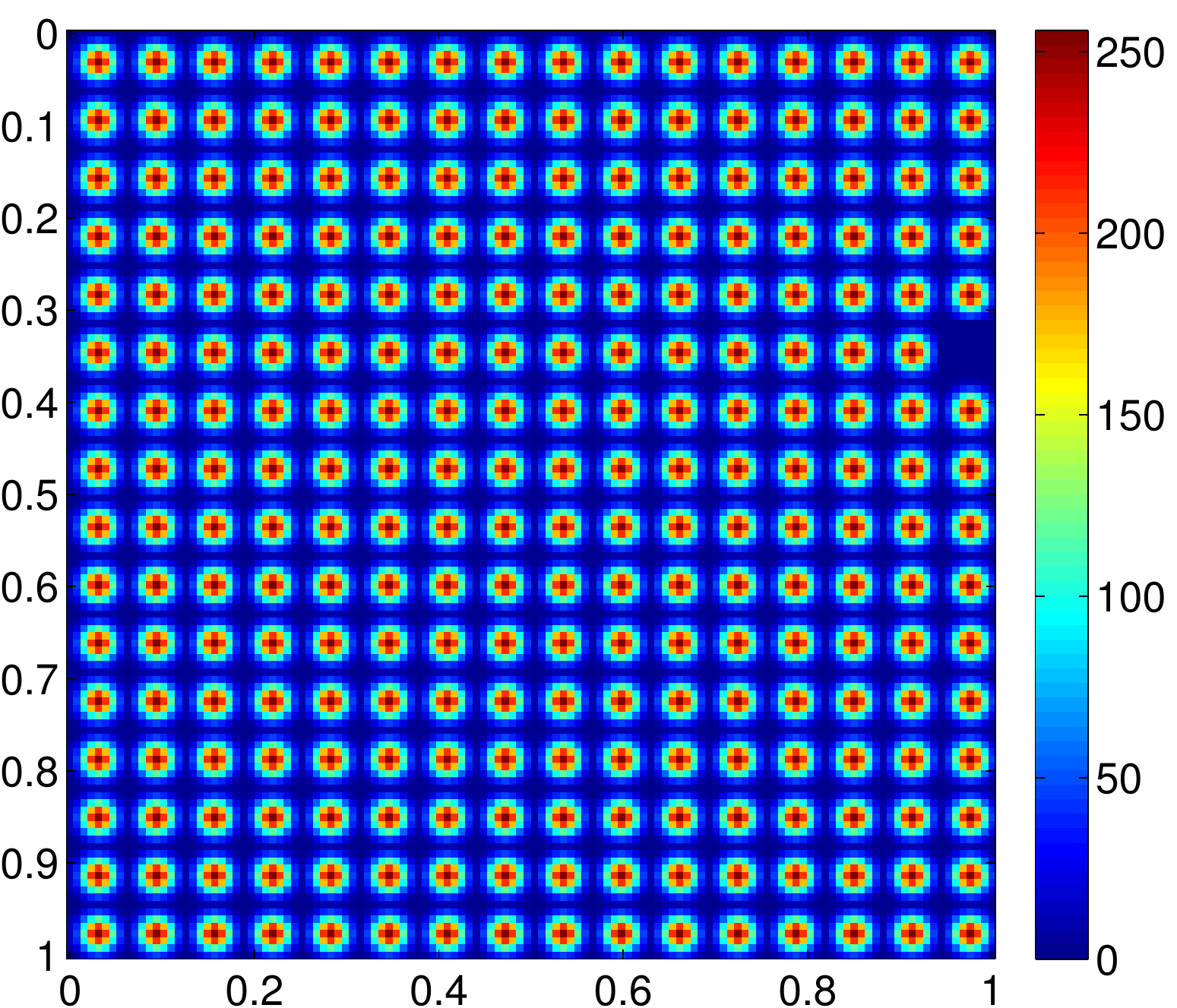} & \includegraphics[height=2.075in]{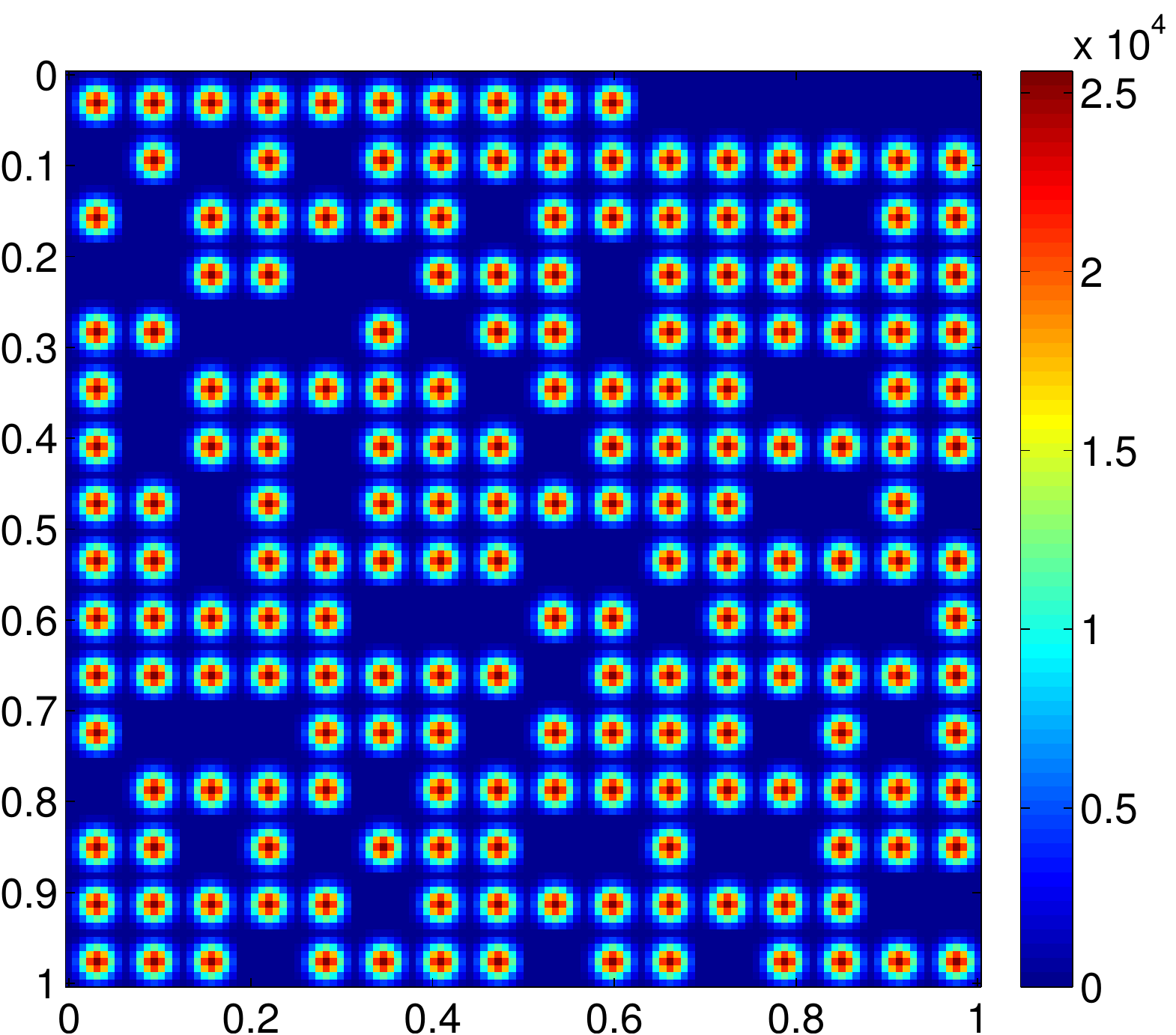} 
\end{tabular}
\end{center}
\caption{Left: the potential energy operator $V(\mathbf{r})$ with $\ell=16$ in Test $2$. Right: the potential energy operator $V(\mathbf{r})$ with $\ell=16$ in Test $3$.}
\label{fig:pot2}
\end{figure}

In the second test, the Hamiltonian matrix $H$ is a discrete
representation of the Hamiltonian operator in a two-dimensional
Kohn-Sham equation similar to the one in \eqref{eqn:ex0}, but the
potential energy operator is more prominent and has a local defect
representing a vacancy of the lattice, i.e., there is a random vacant
cell in $[0,\ell)^2$. The potential energy operator $V(\mathbf{r})$ is
constructed by randomly covering a Gaussian well in $V_0(\mathbf{r})$
with a zero patch (see Figure~\ref{fig:pot2} left panel for an
example).
     
The Hamiltonian matrix $H$ is not dominated by the kinetic matrix
$-\frac{\Delta}{2}$ and hence their spectra are different.  Therefore,
the performance of the regular preconditioners TPA and gTPA in Test
$2$ might not be as good as in Test $1$. For a similar problem size
$n$, the number of iterations in the OMM is significantly larger than
the one in Test $1$. The measurement $d(X,X_0)$ is as large as $1e-4$
meaning that the preconditioned OMM with TPA and gTPA becomes less
accurate in revealing the true eigenspace. Since the number of vacant
cells if fixed, the influence of the local defect becomes weaker as
the number of cells $\ell$ (per dimension) becomes larger. Hence, the
Hamiltonian $H$ becomes closer to the kinetic operator and the number
of iteration decreases. This also shows that the performance of TPA
and gTPA is better for Hamiltonian matrices close to the kinetic part.
     
\begin{figure}
\begin{center}
   \begin{tabular}{cc}
        \includegraphics[height=2in]{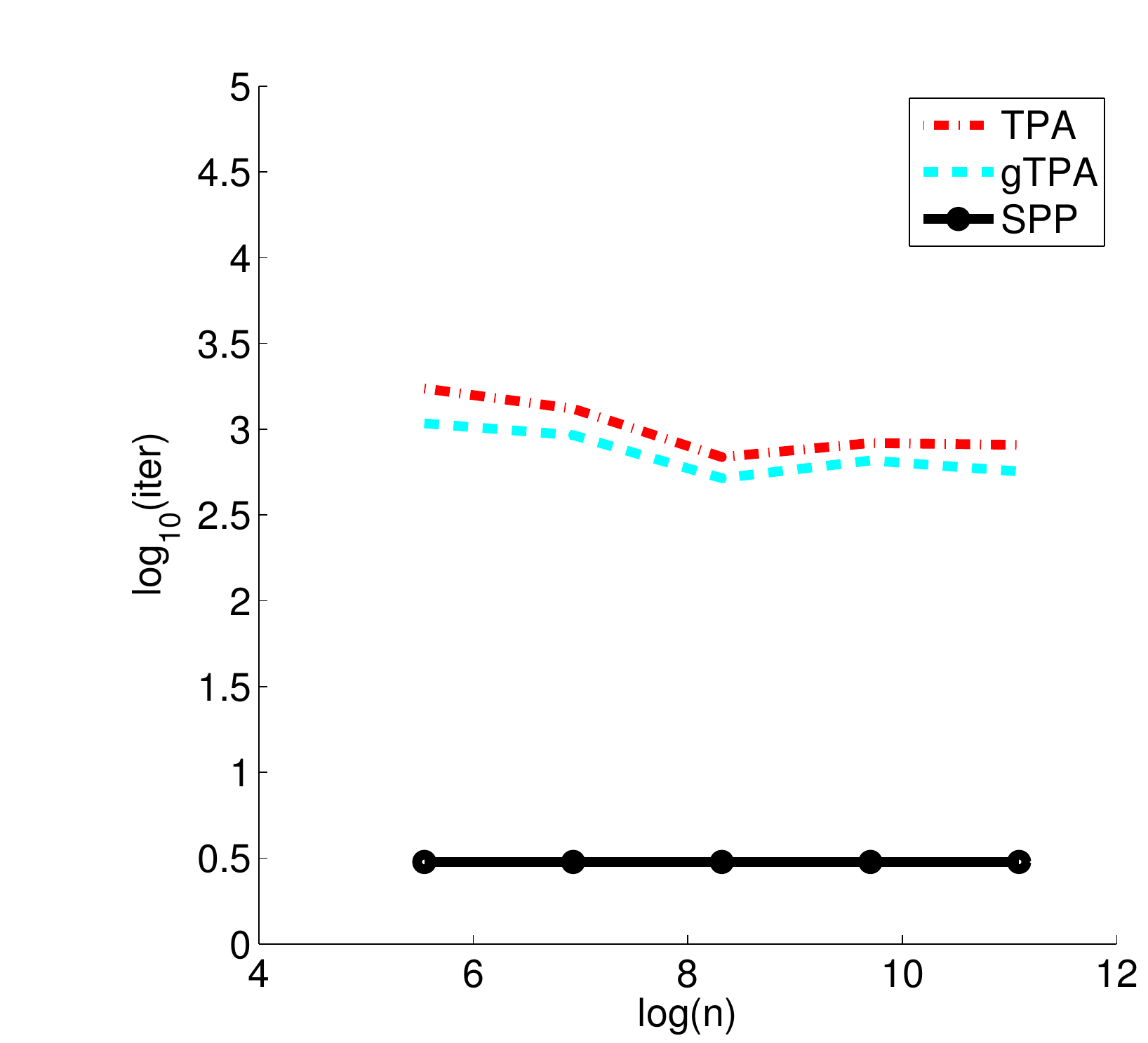} &\includegraphics[height=2in]{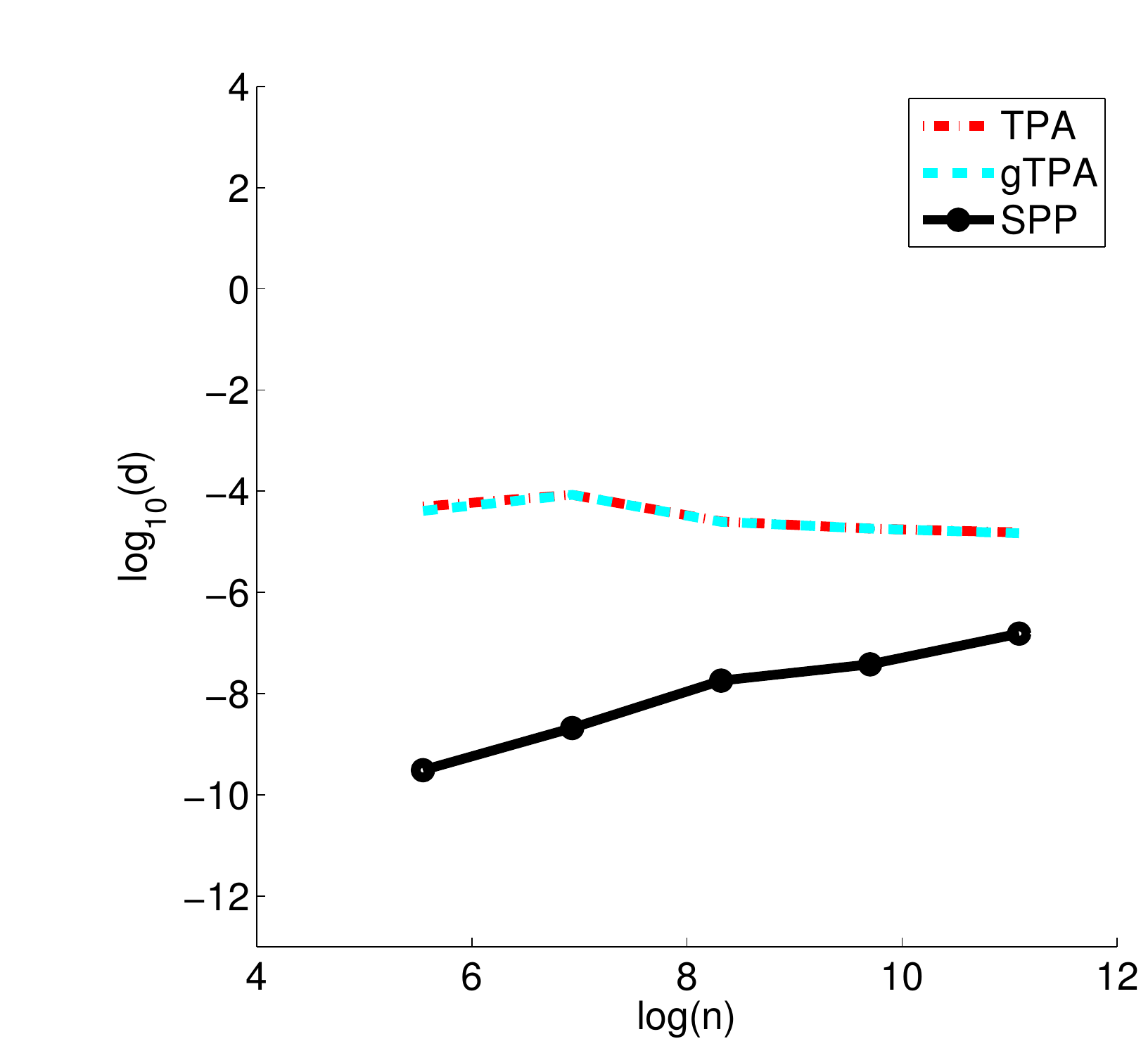}  \\
        (a) & (b) \\
         \includegraphics[height=2in]{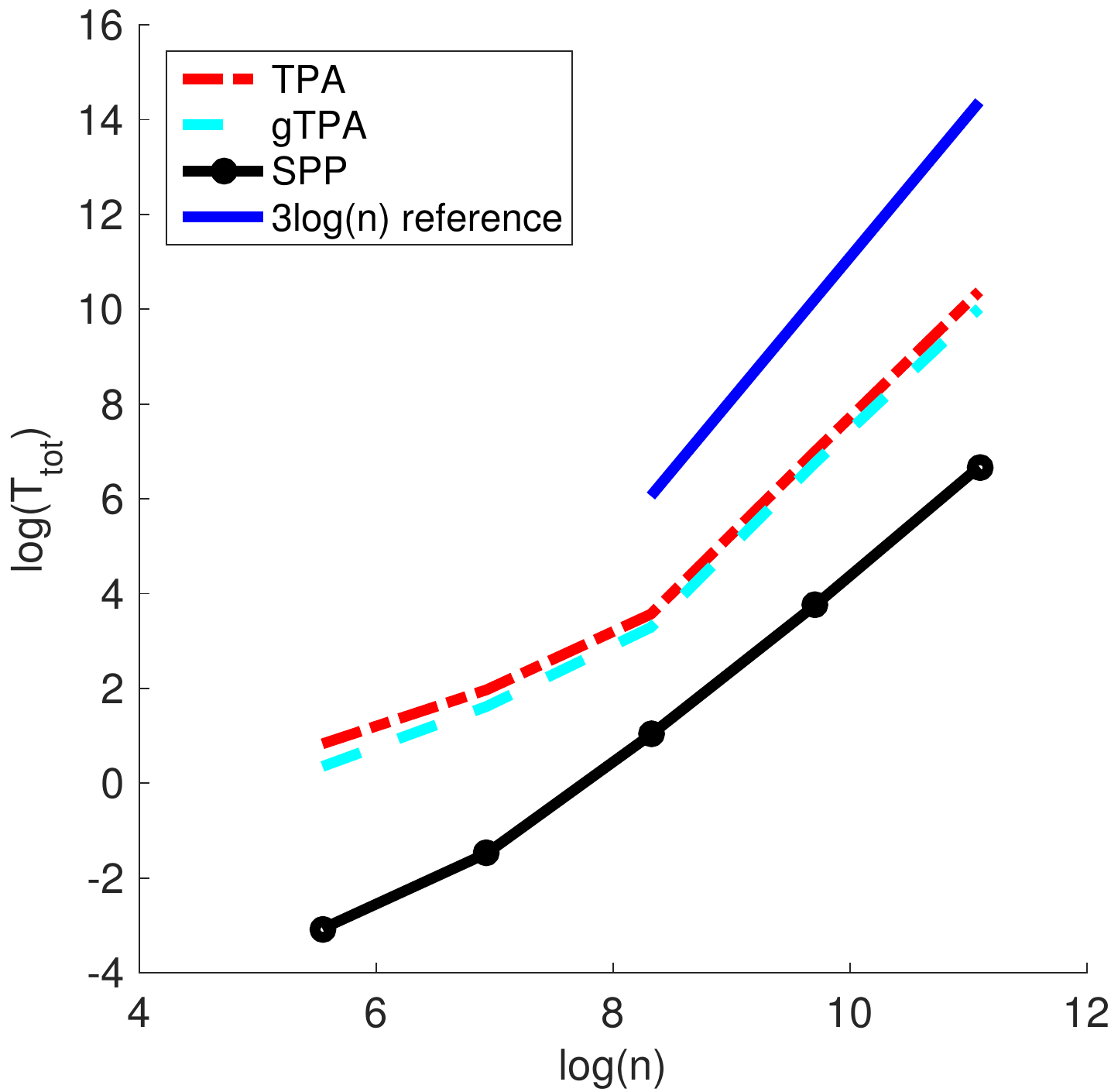} &\includegraphics[height=2in]{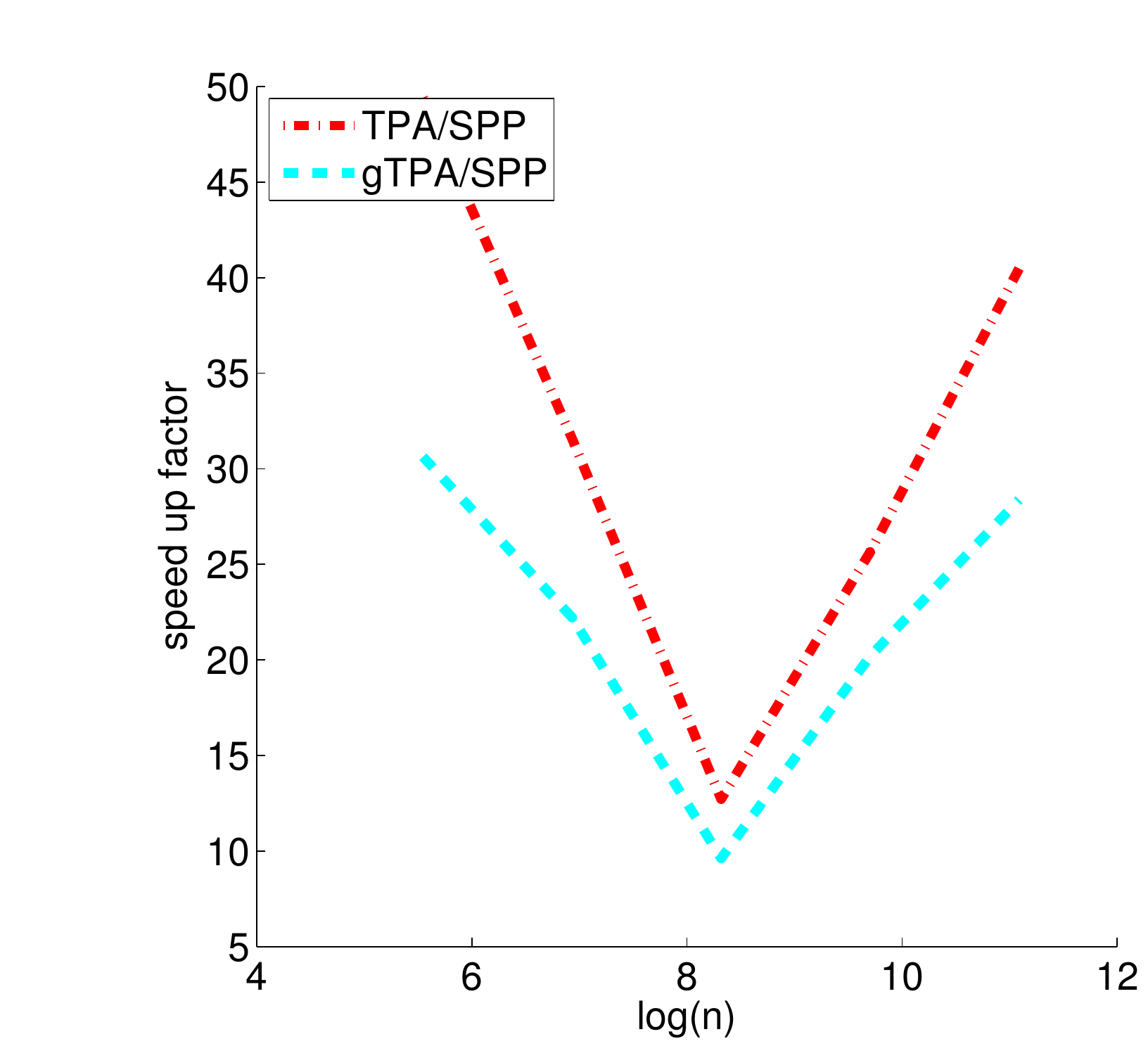} \\
\end{tabular}
\end{center}
\caption{Numerical results in Test $2$. (a) the number of iterations in the preconditioned OMM. (b) the measurement $d(X,X_0)$. (c) the total running time $T_{\text{tot}}$. (d) the speedup factor of the PP method compared with the TPA and gTPA method.}
\label{fig:set2}
\end{figure}

In Test $2$, since the potential energy is prominent and the Hamiltonian operator behaves far from the kinetic operator, the GMRES with an inverse shifted kinetic energy operator is not efficient for solving linear systems like
\[
(H-z_kI)x=b.
\]
The GMRES and the recently developed sparsifying preconditioner \cite{Ying:spspc,Ying:spspd,LuYing:nonSch} are applied to solve the above system in Test $2$ and $3$. The preconditioner via an approximate Fermi operator projection for the OMM is hence denoted as SPP. As shown in Table \ref{tab:2} and Figure \ref{fig:set2}, the performance of the SPP method is better than the TPA and gTPA methods: the number of iteration in the OMM is a small number of $\Or(1)$; the preconditioned OMM with SPP is able to provide an eigenspace closer to the ground truth; the SPP method is much faster and the speedup factor tends to increase with $n$.

\begin{table}[htp]
\centering
\begin{tabular}{rccccccc}
\toprule
    & $(\ell,n)$ & $\mathsf{cond}$ & $\mathsf{iter}$ & $T_{\text{st}}(sec)$
                        &  $T_{\text{omm}}(sec)$   & $T_{\text{tot}}(sec)$ & $d$ \\
\toprule
TPA & (2,256) & 3.7e+03 & 1.7e+03 & 5.514e-04 & 2.295e+00 & 2.296e+00 & 5.0e-05 \\ 
gTPA & (2,256) & 3.7e+03 & 1.1e+03 & 5.518e-04 & 1.423e+00 & 1.423e+00 & 4.1e-05 \\ 
SPP & (2,256) & 3.7e+03 & 3.0e+00 & 4.175e-02 & 4.708e-03 & 4.646e-02 & 3.1e-10 \\ 
\toprule
TPA & (4,1024) & 2.4e+05 & 1.3e+03 & 2.652e-04 & 7.180e+00 & 7.180e+00 & 8.5e-05 \\ 
gTPA & (4,1024) & 2.4e+05 & 9.2e+02 & 3.680e-04 & 5.052e+00 & 5.053e+00 & 8.5e-05 \\ 
SPP & (4,1024) & 2.4e+05 & 3.0e+00 & 2.074e-01 & 2.007e-02 & 2.275e-01 & 2.1e-09 \\ 
\toprule
TPA & (8,4096) & 9.4e+05 & 6.9e+02 & 1.098e-03 & 3.542e+01 & 3.542e+01 & 2.5e-05 \\ 
gTPA & (8,4096) & 9.4e+05 & 5.2e+02 & 1.489e-03 & 2.680e+01 & 2.680e+01 & 2.5e-05 \\ 
SPP & (8,4096) & 9.4e+05 & 3.0e+00 & 2.590e+00 & 1.965e-01 & 2.787e+00 & 1.8e-08 \\ 
\toprule
TPA & (16,16384) & 7.4e+05 & 8.3e+02 & 1.432e-02 & 1.096e+03 & 1.096e+03 & 1.8e-05 \\ 
gTPA & (16,16384) & 7.4e+05 & 6.6e+02 & 1.453e-02 & 8.630e+02 & 8.630e+02 & 1.8e-05 \\ 
SPP & (16,16384) & 7.4e+05 & 3.0e+00 & 3.774e+01 & 5.000e+00 & 4.274e+01 & 3.8e-08 \\ 
\toprule
TPA & (32,65536) & 3.1e+06 & 8.1e+02 & 1.754e-01 & 3.211e+04 & 3.211e+04 & 1.5e-05 \\ 
gTPA & (32,65536) & 3.1e+06 & 5.7e+02 & 2.054e-01 & 2.251e+04 & 2.251e+04 & 1.5e-05 \\ 
SPP & (32,65536) & 3.1e+06 & 3.0e+00 & 6.459e+02 & 1.473e+02 & 7.932e+02 & 1.5e-07 \\

\bottomrule
\end{tabular}
\caption{Numerical results in Test $2$. }
\label{tab:2}
\end{table}

\subsection{Test $3$}
     
In the last test, the potential energy operator in Test $2$ is multiplied by $100$ and $25\%$ of the cells are randomly covered by zero patches (see Figure \ref{fig:pot2} right panel for an example). Different from Test $2$ that keeps the same number of vacant cells and varies the problem size, the number of vacant cells in Test $3$ is proportional to the problem size so that the Hamiltonian matrix $H$ stays far away from the kinetic matrix $-\frac{\Delta}{2}$, and their spectra are very different. Hence, the preconditioned OMM with TPA or gTPA is no longer efficient for all problem sizes. 

\begin{figure}
\begin{center}
   \begin{tabular}{cc}
                 \includegraphics[height=2in]{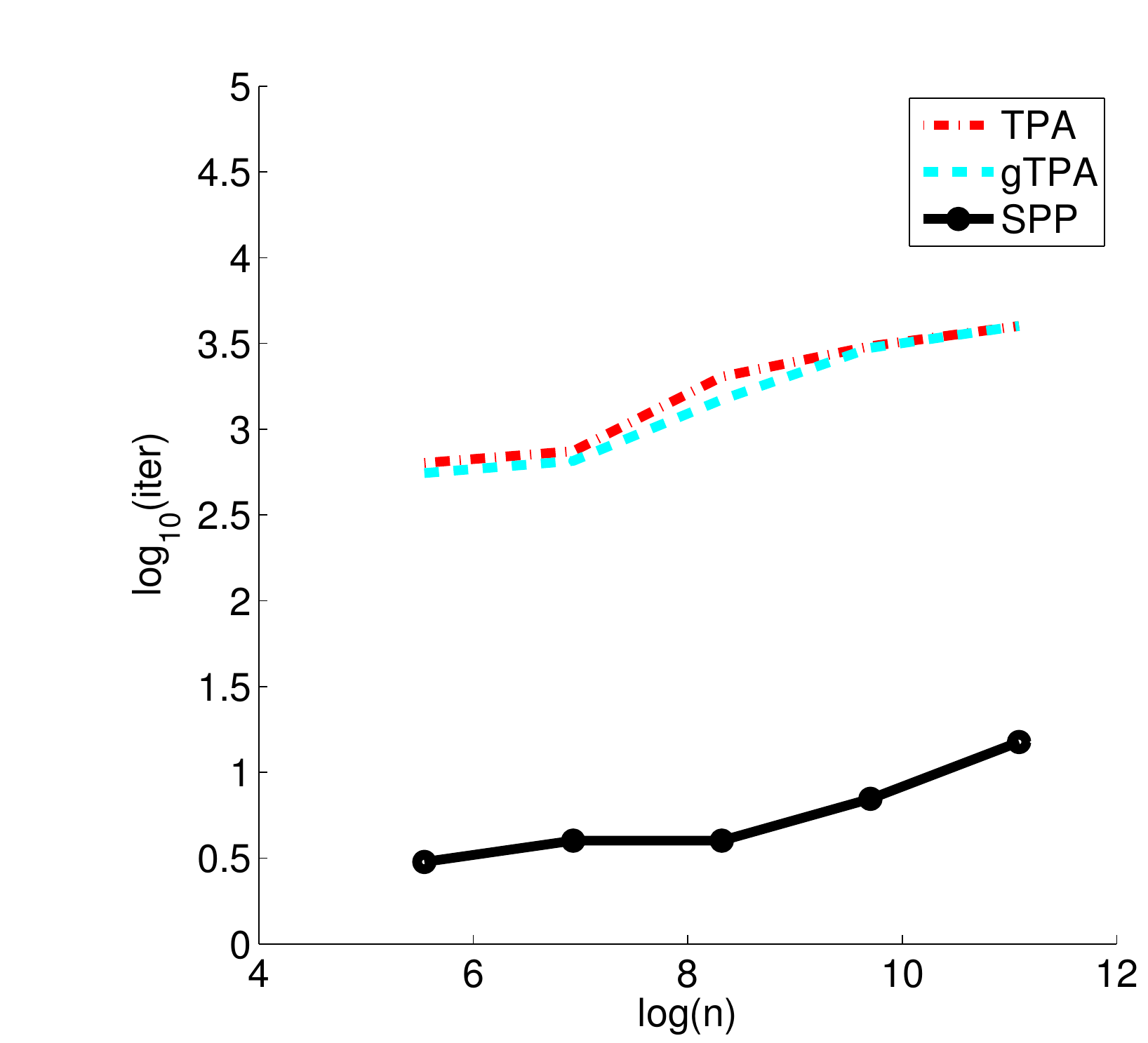} &\includegraphics[height=2in]{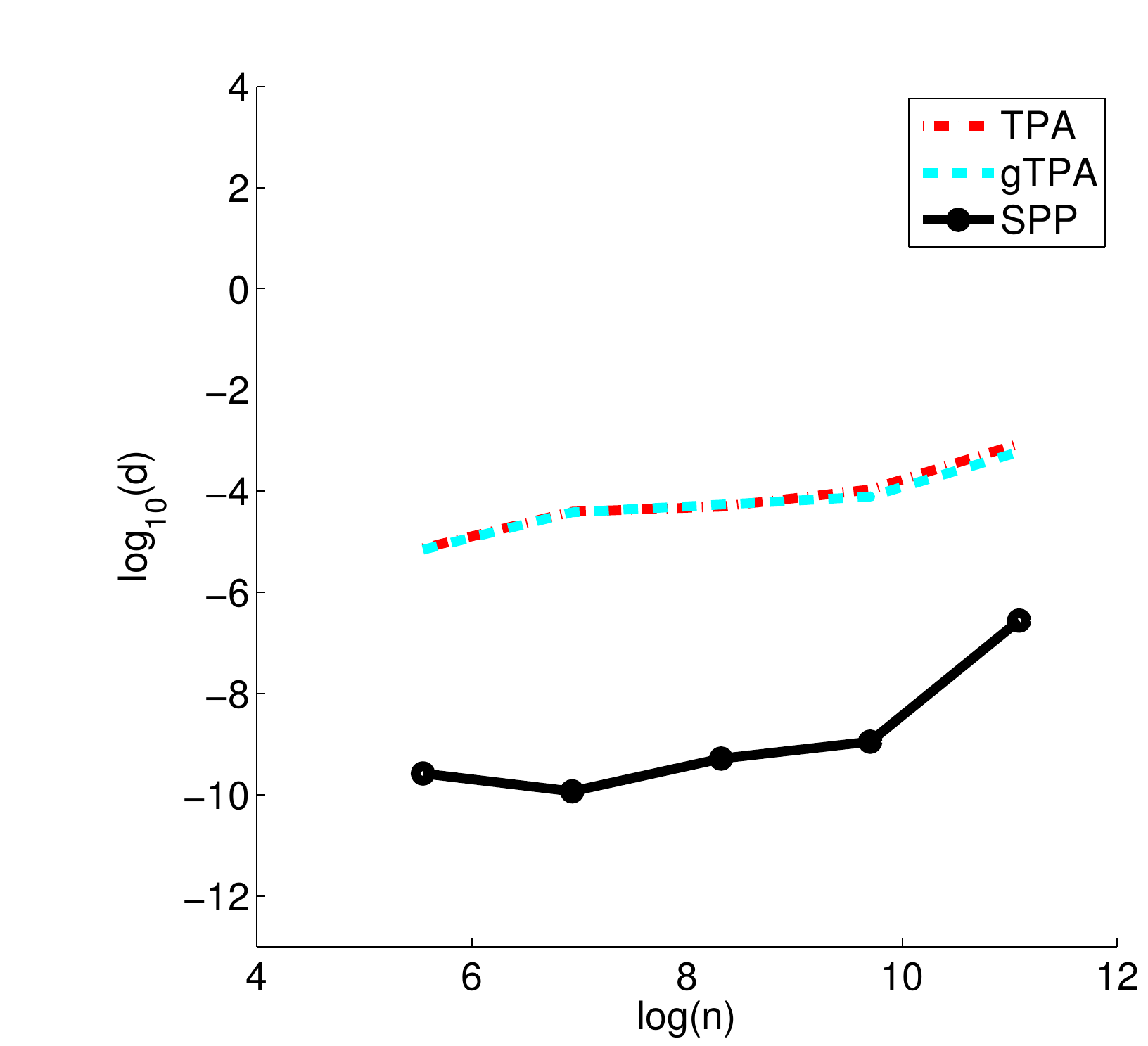}  \\
        (a) & (b) \\
         \includegraphics[height=2in]{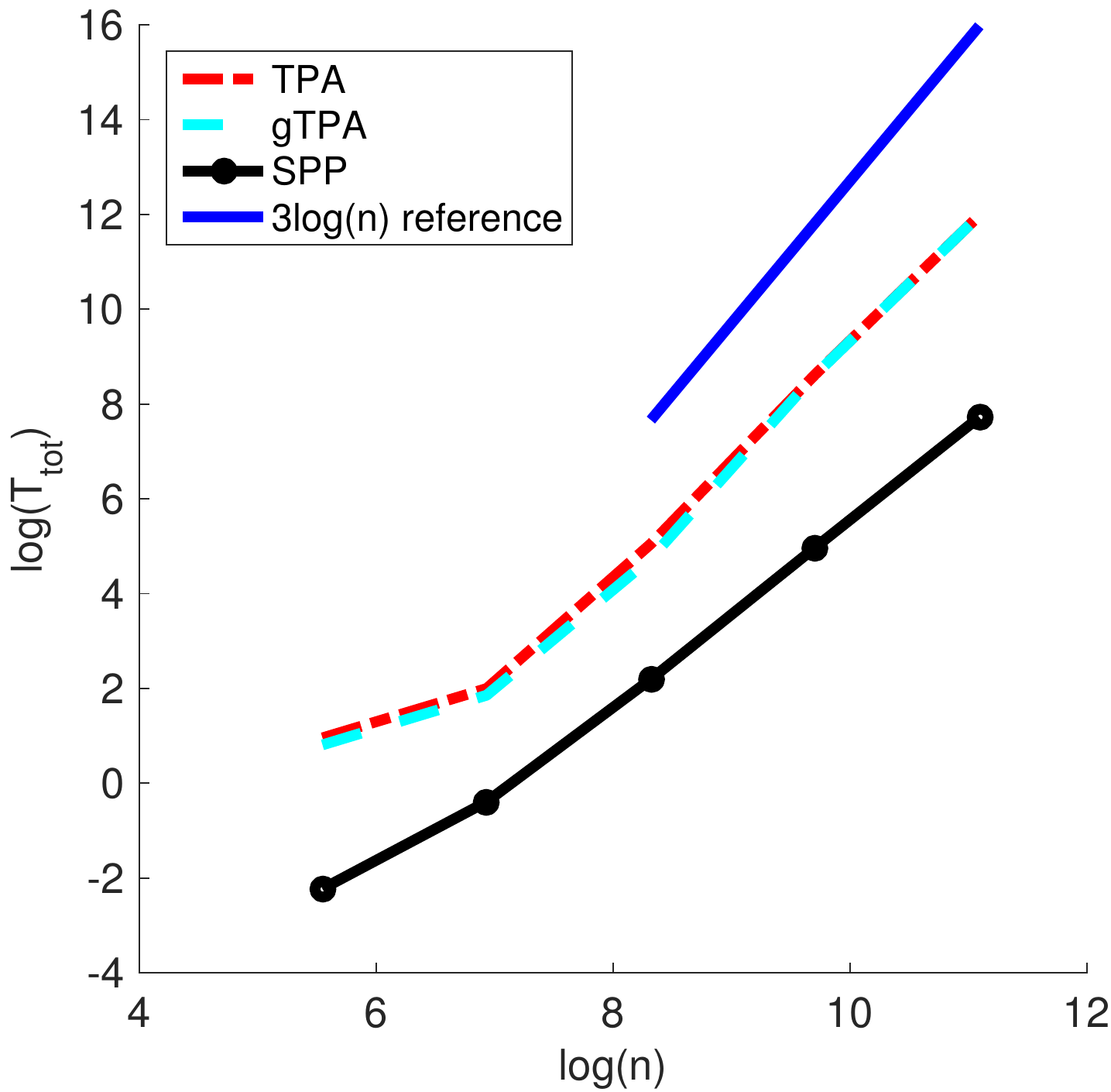}& \includegraphics[height=2in]{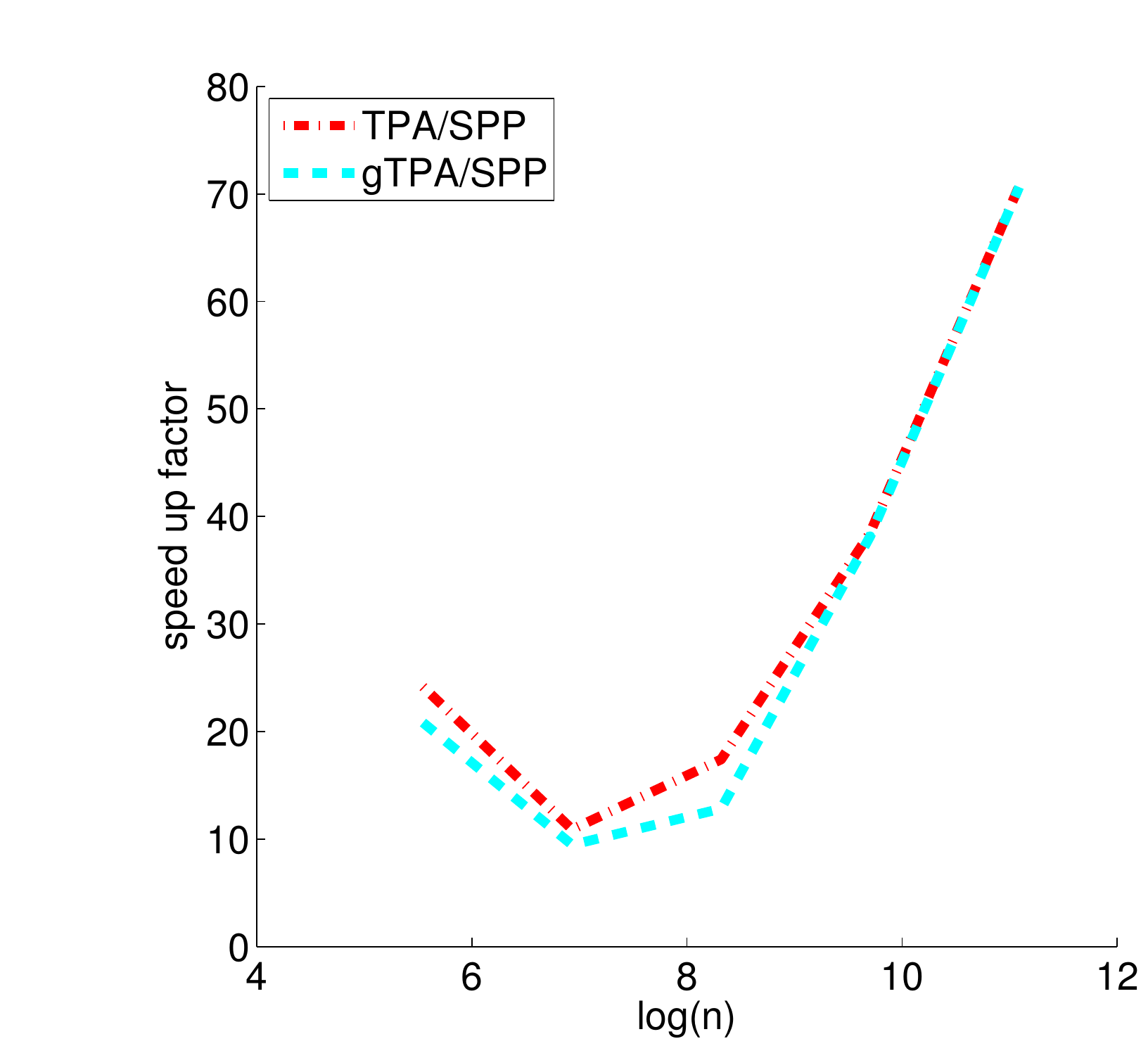} \\
        (c) & (d) 
\end{tabular}
\end{center}
\caption{Numerical results in Test $3$. (a) the number of iterations in the preconditioned OMM. (b) the measurement $d(X,X_0)$. (c) the total running time $T_{\text{tot}}$. (d) the speedup factor of the PP method compared with the TPA and gTPA method.}
\label{fig:set3}
\end{figure}

As shown in Table \ref{tab:3} and Figure \ref{fig:set3}, the preconditioned OMM with TPA or gTPA requires thousands of iterations to converge or cannot converge within the maximum number of iteration, $4000$. Even if the OMM converges, it cannot return useful estimation of the eigenspace, since the measurement $d(X,X_0)$ is too large.

 In contrast, the OMM with SPP is still able to provide reasonably good estimation of the eigenspace within $\Or(1)$ iterations. Although the iteration number increases with the system size $n$, it remains quite small even for large-scale problems.

\begin{table}[htp]
\centering
\begin{tabular}{rccccccc}
\toprule
    & $(\ell,n)$ & $\mathsf{cond}$ & $\mathsf{iter}$ & $T_{\text{st}}(sec)$
                        &  $T_{\text{omm}}(sec)$   & $T_{\text{tot}}(sec)$ & $d$ \\
\toprule
TPA & (2,256) & 6.8e+02 & 6.4e+02 & 1.098e-03 & 2.602e+00 & 2.603e+00 & 7.4e-06 \\ 
gTPA &  (2,256) & 6.8e+02 & 5.5e+02 & 9.694e-04 & 2.244e+00 & 2.245e+00 & 7.0e-06 \\ 
SPP & (2,256) & 6.8e+02 & 3.0e+00 & 9.940e-02 & 8.341e-03 & 1.077e-01 & 2.7e-10 \\ 
\toprule
TPA & (4,1024) & 2.5e+03 & 7.5e+02 & 6.446e-04 & 7.374e+00 & 7.374e+00 & 3.9e-05 \\ 
gTPA & (4,1024) & 2.5e+03 & 6.5e+02 & 4.654e-04 & 6.401e+00 & 6.402e+00 & 3.8e-05 \\ 
SPP & (4,1024) & 2.5e+03 & 4.0e+00 & 6.153e-01 & 6.153e-02 & 6.768e-01 & 1.2e-10 \\ 
\toprule
TPA & (8,4096) & 3.7e+03 & 2.0e+03 & 2.367e-03 & 1.562e+02 & 1.562e+02 & 4.9e-05 \\ 
gTPA & (8,4096) & 3.7e+03 & 1.5e+03 & 2.307e-03 & 1.151e+02 & 1.151e+02 & 5.4e-05 \\ 
SPP & (8,4096) & 3.7e+03 & 4.0e+00 & 8.540e+00 & 4.380e-01 & 8.978e+00 & 5.2e-10 \\ 
\toprule
TPA & (16,16384) & 2.4e+04 & 3.0e+03 & 2.196e-02 & 5.557e+03 & 5.557e+03 & 1.1e-04 \\ 
gTPA & (16,16384) & 2.4e+04 & 3.0e+03 & 2.476e-02 & 5.509e+03 & 5.509e+03 & 7.8e-05 \\ 
SPP & (16,16384) & 2.4e+04 & 7.0e+00 & 1.288e+02 & 1.563e+01 & 1.444e+02 & 1.1e-09 \\ 
\toprule
TPA & (32,65536) & 1.4e+05 & - & 2.016e-01 & 1.595e+05 & 1.595e+05 & 8.7e-04  \\ 
gTPA & (32,65536) & 1.4e+05 & - & 2.308e-01 & 1.593e+05 & 1.593e+05 & 6.0e-04 \\ 
SPP & (32,65536) & 1.4e+05 & 1.5e+01 & 1.578e+03 & 6.779e+02 & 2.256e+03 &2.8e-07 \\ 
 
\bottomrule
\end{tabular}
\caption{Numerical results in Test $3$. ``$-$" means exceeding the maximum iteration number $4000$.}
\label{tab:3}
\end{table}

\section{Conclusion and Remarks} 
\label{sec:conclusion}
This paper presents a novel preconditioner for the orbital minimization method (OMM) from planewave discretization. Once constructed, the application of the preconditioner is very efficient and the preconditioned OMM converges in $\Or(1)$ iterations. Based on the approximate Fermi operator projection, this preconditioner can be constructed efficiently via iterative matrix solvers like GMRES with the newly developed sparsifying preconditioner. The speedup factor of the running time compared with existing preconditioned OMMs is as large as hundreds of times and the speedup factor tends to increase with the problem size. Numerical experiments also show that the new preconditioned OMM is able to provide more accurate solutions than the popular TPA preconditioner and, as a result, might reduce the number of iterations in the self-consistent field iteration.

This preconditioner based on spectral projection can be considered for other discretization schemes using localized basis functions, e.g., numerical atomic orbitals and wavelets that lead to sparse Hamiltonians. In these cases, sparse direct solvers or preconditioned iterative solvers could be applied to construct the approximate Fermi operator projection via the pole expansion. This would be worth exploring as future directions.

It would be also of interest to implement our
new algorithm into the recently developed parallel library for the OMM
framework, \textsf{libomm} \cite{Corsetti:14} and incorporate into
existing electronic structure software packages.  Parallelism of this
new preconditioned OMM is straightforward since the main routines of
this algorithm have parallel analog in existing high performance
computing packages. Several versions of parallel FFT can be found in
\cites{PippigPotts:13,Poulson:14}; the computational cost per
processor is less in the former one, while the latter one has a better
scalability and the number of processor can be $\Or(n)$; a recent
article \cite{LiYang:16} reduces the prefactor of the algorithm in
\cite{Poulson:14} to an optimal one in the butterfly scheme while
keeping the same scalability. The pole expansion can be embarrassingly
parallelized. The last main routine, if applied, is the scalable QR
factorization used for precomputing and storing the preconditioner in
a data-sparse format. This has been implemented in Elemental
\cite{Poulson:2013}. For plane wave discretization, a recent article
\cite{Levitt:15} has shown that the constrained minimization model in
\eqref{eq:tracemin} with a Chebyshev filter \cite{Zhou2006}
outperforms the locally optimal block preconditioned conjugate
gradient (LOBPCG) algorithm \cite{Knyazev:2001}, implemented in ABINIT
\cite{Bottin2008329}, in terms of scalability, even though a direct
diagonlization of a matrix of size $N\times N$ is required in the
Chebyshev filter in each iteration.  Since the new preconditioned OMM
requires no direct diagonlization, it is promising to lead to a more
scalable approach.

\bibliographystyle{abbrv}
\bibliography{nonSch}

\end{document}